\documentclass[11pt, twoside]{article}

\usepackage[english]{babel}

\usepackage{bbm}
\usepackage{amssymb}
\usepackage{amsfonts}
\usepackage{amsmath}
\usepackage{amsthm}
\usepackage{color}
\usepackage{mathrsfs}
\usepackage{txfonts}
\usepackage{bbm}
\usepackage{enumerate}
\usepackage{anysize}
\usepackage{indentfirst}

\usepackage{latexsym}

\usepackage[colorlinks=true,
linkcolor=blue,
citecolor=red,
urlcolor=magenta,
backref=page
]{hyperref}

\allowdisplaybreaks

\newtheorem{theorem}{Theorem}[section]
\newtheorem{lemma}[theorem]{Lemma}

\newtheorem{proposition}[theorem]{Proposition}

\theoremstyle{definition}
\newtheorem{remark}[theorem]{Remark}
\newtheorem{definition}[theorem]{Definition}

\newcounter{assum}

\renewcommand{\appendix}{\par
	\setcounter{section}{0}
	\setcounter{subsection}{0}
	\setcounter{subsubsection}{0}
	\gdef\thesection{\@Alph\c@section}
	\gdef\thesubsection{\@Alph\c@section.\@arabic\c@subsection}
	\gdef\theHsection{\@Alph\c@section.}
	\gdef\theHsubsection{\@Alph\c@section.\@arabic\c@subsection}
	\csname appendixmore\endcsname
}

\numberwithin{equation}{section}

\begin{document}
\title{\bf\Large On VMO-Type Characterizations of
Ball Banach Sobolev Spaces and Independence of Integrable Index of Special
John--Nirenberg--Campanato Spaces
\footnotetext{\hspace{-0.35cm} 2020
\emph{Mathematics Subject Classification}. Primary 46E35;
Secondary 42B25, 26D10, 42B35, 47A30.
\endgraf \emph{Key words and phrases}.
Sobolev space, ball Banach function space,
special John--Nirenberg--Campanato space,
BBM formulas, ball average operator.
\endgraf
}}

\author{Zhexu Bai, Hao Peng,
Dachun Yang,
Wen Yuan and Yangyang Zhang}
\maketitle

\vspace{-0.8cm}

\begin{center}
\begin{minipage}{13cm}
{\small {\bf Abstract}\quad
Let $k\in\mathbb N$, $X$ be a ball Banach function space, and
$W^{k,X}$ denote the corresponding ball Banach Sobolev space.
In this article, we solve six questions in [Math. Z. 311 (2025)].
Specifically, we first prove that the constants in its VMO-type
characterization of $W^{k,X}$
are sharp uniformly over all ball Banach function spaces
satisfying suitable assumptions.
Next, we show that the BBM-type formula and the BSVY-type formula
involving minimizing polynomials or ball averages fail in general.
In contrast, we establish higher-order BBM-type formulas
involving Taylor polynomials and further obtain a
Taylor-type characterization of $W^{k,X}$ under suitable assumptions.
Moreover, we construct counterexamples to show that the
characterization of $W^{k,X}$ via higher-order
minimizing polynomials fails in general.
Finally, we provide a necessary and sufficient condition
for special John--Nirenberg--Campanato spaces defined via congruent
cubes to be independent of their integrability exponent $q$.}
\end{minipage}
\end{center}


\tableofcontents

\vspace{0.2cm}

\section{Introduction}\label{in}

Throughout this article, we work in $\mathbb{R}^n$ and,
unless otherwise specified, we take $\mathbb{R}^n$ as the default
underlying space.

The Bourgain--Brezis--Mironescu (BBM) \emph{formula}
was first established by Bourgain et al. \cite{bbm01} to
characterize Sobolev spaces on smooth bounded domains in $\mathbb{R}^n$.
Brezis \cite{b02} subsequently extended this formula to $\mathbb{R}^n$:
for any $p\in[1,\infty)$ and $f\in\dot{W}^{1,p}\cap L^p$
(see Definition \ref{def-sobolev}
for the definition of $\dot{W}^{1,p}$),
\begin{align*}
\lim_{s\to1^-}(1-s)
\int_{\mathbb{R}^n}\int_{\mathbb{R}^n}
\frac{|f(x)-f(y)|^p}{|x-y|^{n+sp}}
\,dy\,dx=\frac{K_{p,n}}{p}
\int_{\mathbb{R}^n}|\nabla f(x)|^p\,dx,
\end{align*}
where $s\to 1^-$ means $s\in(0,1)$ and $s\to 1$
and where $K_{p,n}$ is a positive constant depending only on
$p$ and $n$. From then on, the BBM formula has received a lot of attention
and found wide applications in characterizing Sobolev-type spaces;
see, for example, \cite{bn06,bn16,bn18,dlyyz2024,dltyy24,dltyy25,dm23}.
As a weak-type counterpart of the BBM formula,
Brezis et al. \cite{bsvy24} established the following
Brezis--Seeger--Van Schaftingen--Yung (BSVY)
\emph{formula}:
for any $p\in[1,\infty)$ and $f\in C_{\mathrm{c}}^{\infty}$,
\begin{align*}
\lim_{\lambda\to\infty}
\lambda^p\left|\left\{(x,y)\in\mathbb{R}^{2n}:
\frac{|f(x)-f(y)|}{|x-y|^{1+\frac{n}{p}}}>\lambda\right\}
\right|=\frac{C_{p,n}}{n}
\left\| \ \left|\nabla f\right|\ \right\|_{L^p}^p,
\end{align*}
where $C_{p,n}$ is a positive constant depending only on $p$ and $n$.

Rather than varying the smoothness exponent as in the BBM formula,
using local mean oscillations over families of small cubes,
Bourgain et al. \cite{bbm15} defined a new BMO-type space $\mathcal{B}$
on the unit cube $[0,1)^n$ and showed that any
integer-valued function in the VMO-type
subspace of $\mathcal{B}$ should be a constant function.
We refer to \cite{bn95} for more properties concerning the VMO space on $[0,1)^n$.
Based on the work of Bourgain et al. \cite{bbm15}, Fusco et al. \cite{FMS(2018)}
introduced a BMO-type
seminorm and studied its relationship with
the Sobolev space $\dot{W}^{1,p}$ for any $p\in(1,\infty)$.
To be precise, Fusco et al. \cite{FMS(2018)} proved that there exists
a positive constant $\gamma_{n,p}$ such that, for any $f\in\dot{W}^{1,p}$,
\begin{align*}
\lim_{\varepsilon\to 0^+}
\sup_{\{Q_j\}_{j\in\mathcal{G}_{\varepsilon}}}
\varepsilon^{n-p}
\sum_{j\in\mathcal{G}_{\varepsilon}}
\fint_{Q_{j}}\left|f(x)-f_{Q_j}\right|^{p}\,dx
=\gamma_{n,p}\|f\|_{\dot{W}^{1,p}}^p,
\end{align*}
where the supremum is taken over all families
$\{Q_j\}_{j\in\mathcal{G}_\varepsilon}$
of pairwise disjoint open cubes
with the same edge length $\varepsilon$ but arbitrary orientation.
Here, and thereafter, the \emph{notation
$\varepsilon\to0^+$} means that $\varepsilon\in(0,1)$
and $\varepsilon\to0$. For any $f\in L_{\mathrm{loc}}^1$
(the set of all locally integrable functions on $\mathbb{R}^n$)
and any cube or any ball $Q\subset\mathbb{R}^n$,
we define the \emph{average of $f$ over $Q$} by
\begin{align*}
f_Q:=\fint_{Q}f(x)\,dx:=\frac{1}{|Q|}\int_{Q}f(x)\,dx.
\end{align*}
The above formula has been further studied and developed;
see, for example, \cite{df20,ffgs20,fms16,gs23}.
In particular, Jia et al. \cite{jyy25} established the
VMO-type characterization of ball Banach Sobolev spaces.
Precisely speaking, let $q\in[1,\infty)$,
$k$ be a positive integer, $X$ be a ball Banach function space
(see Definition \ref{def-ball-Banach}),
and $\dot{W}^{k,X}$ be the ball Banach Sobolev space
(see Definition \ref{def-sobolev}).
Under some mild assumptions on $X$,
Jia et al. \cite{jyy25} found two positive
constants $\gamma(n,q,k)$ and $\Gamma(n,q,k)$ [see \eqref{gandGamma}],
depending only on $n$, $q$, and $k$,
such that, for any $f\in\dot{W}^{k,X}$,
\begin{align}\label{equ-polynomial}
\gamma{(n,q,k)}\left\|\,\left|\nabla^kf\right|\,
\right\|_X\leq
\lim_{r\rightarrow 0^+}\left\|\frac{1}{r^k}
\left[\fint_{B(\cdot,r)}
\left|f(y)-P^{(k-1)}_{B(\cdot,r)}(f)(y)
\right|^q\,dy\right]^{\frac{1}{q}}\right\|_X
\leq\Gamma{(n,q,k)}\left\|\,\left|\nabla^kf\right|
\,\right\|_X,
\end{align}
where, for any $x\in\mathbb{R}^n$ and $r\in(0,\infty)$,
$B(x,r):=\{y\in\mathbb{R}^n:|y-x|<r\}$
and $P_{B(x,r)}^{(k-1)}(f)$ denotes the $(k-1)$-th
minimizing polynomial of $f$ on $B(x,r)$ (see Section \ref{s3}).
In addition, they also showed that, for any $f\in W^{2k,X}$,
\begin{align}\label{equ-ball-average}
\lim_{r\to 0^+}
\left\|\frac{1}{r^{2k}}
\left[\fint_{B(\cdot,r)}
\left|f(y)-B_{k,r}f(y)\right|^q\,dy
\right]^{\frac{1}{q}}\right\|_X
=a_k\left\|\Delta^k f\right\|_X,
\end{align}
where $B_{k,r}$ is a generalization of the ball average operator (see Section \ref{s8}),
$a_k$ is a positive constant depending only on $k$ and $n$,
and $\Delta:=\sum_{j=1}^n\frac{\partial^2}{\partial x_j^2}$ is the \emph{Laplace operator}.
In the same article, Jia et al. asked six questions.
Specifically, in establishing \eqref{equ-polynomial} and \eqref{equ-ball-average},
they left five questions, respectively, in
\cite[Remarks 3.12(ii), 3.14, 3.19(i) and (ii), and 4.10]{jyy25}.
Moreover, they applied \eqref{equ-polynomial} to study
the dependence of the special
John--Nirenberg--Campanato spaces via congruent cubes $JN^{\rm{con}}_{(p,q,k)_\alpha}$
(see Definition \ref{def-john}) on the integrability exponent $q$
in \cite[Subsection 5.6]{jyy25}
and asked a related question in \cite[Remark 5.20]{jyy25}.

In this article, we answer all the above six questions.
Specifically, we establish the sharpness of the constants
$\gamma(n,q,k)$ and $\Gamma(n,q,k)$
in the VMO-type characterizations [see \eqref{equ-polynomial}],
which gives an affirmative answer to the question in
\cite[Remark 3.12(ii)]{jyy25}. We
show that the BBM-type and
BSVY-type formulas involving minimizing polynomials or ball
averages fail in general (see Theorems \ref{wenti1fouding},
\ref{ruofanshu1}, \ref{pro-laplace}, and
\ref{pro-weaknorm}), which gives negative answers to the questions
in \cite[Remarks 3.14 and 4.10]{jyy25}.
We establish the BBM-type
formula involving Taylor polynomials (see Theorem \ref{taile}) and
use it to obtain a new VMO-type characterization
of ball Banach Sobolev spaces by Taylor polynomials
(see Theorem \ref{tailek}), which gives an affirmative
answer to the question in \cite[Remark 3.19(i)]{jyy25}.
We also disprove the proposed characterization involving higher-order
minimizing polynomials (see Theorem \ref{the-4}), which gives a negative
answer to the question in \cite[Remark 3.19(ii)]{jyy25}.
Finally, we prove that the special
John--Nirenberg--Campanato spaces via congruent cubes $JN^{\rm{con}}_{(p,q,k)_\alpha}$
are independent of the integrability exponent $q$ if and
only if $\alpha\in (0,\infty)$ (see Theorem \ref{q-wug}), which answers the question
in \cite[Remark 5.20]{jyy25}.

To present our main results, we first recall some basic concepts.
For any $k\in\mathbb{Z}_+$, let the \emph{notation $C^k$} denote the set of all
$k$-times continuously differentiable functions on $\mathbb{R}^n$,
and let the \emph{notation $C^\infty$} denote the set of all infinitely differentiable
functions on $\mathbb{R}^n$.
Let the \emph{notation $C_{\rm{c}}$} denote the set of all continuous functions on
$\mathbb{R}^n$ with compact support, and, for any
$k\in\mathbb{Z}_+\cup\{\infty\}$, let $C_{\rm{c}}^k:=C^k\cap C_{\rm{c}}$.

Our first result concerns the sharpness of the constants
$\gamma(n,q,k)$ and $\Gamma(n,q,k)$ in \eqref{equ-polynomial}.
Jia et al. proved that these constants are sharp in some special cases
(see \cite[Remark 3.12(iii)]{jyy25}).
By constructing suitable weights and
weighted Lebesgue spaces (see Theorem \ref{wt1jd}),
we prove that both constants are sharp uniformly over
all ball Banach function spaces $X$ satisfying the assumptions
in \cite[Theorem 3.3(i)]{jyy25}.
\begin{theorem}\label{zjx}
Let $n,k\in\mathbb{N}$ and $q\in[1,\infty)$.
Then
$$\gamma(n,q,k)=\inf_{X}\inf_{\genfrac{}{}{0pt}{}{f\in C^k,|\nabla^kf|\in C_{\rm{c}}}
{\|\,|\nabla^kf|\,\|_{X}\neq0}}
\frac{\|\mathrm{Osc}^{(q)}_k(\nabla(D^{k-1}f))\|_X}
{\|\,|\nabla^kf|\,\|_X}$$
and
$$\Gamma(n,q,k)=\sup_{X}\sup_{\genfrac{}{}{0pt}{}{f\in C^k,|\nabla^kf|\in C_{\rm{c}}}
{\|\,|\nabla^kf|\,\|_{X}\neq0}}
\frac{\|\mathrm{Osc}^{(q)}_k(\nabla(D^{k-1}f))\|_X}
{\|\,|\nabla^kf|\,\|_X},$$
where $\gamma(n,q,k)$ and $\Gamma(n,q,k)$ are as in \eqref{gandGamma},
$\mathrm{Osc}^{(q)}_k(A)$ and $|\nabla^kf|$ are
as, respectively, in \eqref{OSCD} and \eqref{def-wf},
and  $\nabla g$ and $D^kf$ are as in Section \ref{s3}.
Here the infimum and the supremum
over $X$ are taken over all ball Banach function spaces $X$
such that $X^{\frac{1}{q}}$ (see Definition \ref{tuhua}) is also a ball Banach function space.
\end{theorem}

\begin{remark}
Theorem \ref{zjx} gives an affirmative answer to the question
in \cite[Remark 3.12(ii)]{jyy25}. However,
we point out that the uniform sharpness in Theorem \ref{zjx}
does not imply the sharpness of these constants for
any fixed ball Banach function space $X$ (see Proposition \ref{exp-3.4}).
\end{remark}

We next investigate the BBM-type formulas proposed in \cite[Remarks 3.14 and 4.10]{jyy25}
and their corresponding BSVY-type analogues.
Theorems \ref{wenti1fouding} and \ref{pro-laplace} give negative answers
to the questions, respectively, in \cite[Remarks 3.14 and 4.10]{jyy25}.
Moreover, Theorems \ref{ruofanshu1} and \ref{pro-weaknorm} show that
the corresponding BSVY-type analogues also fail in general.
Specifically, these results are obtained by constructing a family of singular
radial kernels [see \eqref{teshu}]. We first consider
the BBM-type formula involving minimizing polynomials.
\begin{theorem}\label{wenti1fouding}
Let $q\in[1,\infty)$ and $\theta\in(0,\min\{n,q\})$.
For any $r\in(0,\infty)$ and $x\in\mathbb{R}^n$, define
\begin{equation}\label{teshu}
\rho_r(x):=
\frac{\theta\mathbf{1}_{B(\mathbf{0},r)
\setminus\{0\}}(x)}{\sigma(\mathbb{S}^{n-1})
r^\theta|x|^{n-\theta}},
\end{equation}
where $\sigma(\mathbb{S}^{n-1})$ denotes the surface Lebesgue measure
of the unit sphere $\mathbb{S}^{n-1}\subset\mathbb{R}^n$.
Then, for any $f\in C_{\rm{c}}^2$,
\begin{equation}\label{equ-infty}
\liminf_{r\to 0^+}
\left[\int_{\mathbb{R}^n}
\int_{\mathbb{R}^n}
\frac{|f(y)-f_{B(x,r)}|^q}{|x-y|^q}
\rho_r(x-y)\,dy\,dx\right]^{\frac{1}{q}}<\infty
\end{equation}
if and only if $f\equiv 0$.
\end{theorem}
The corresponding BSVY-type analogue involving minimizing polynomials also fails in general
as the next theorem shows.
\begin{theorem}\label{ruofanshu1}
Let $q\in[1,\infty)$ and $\{\rho_r\}_{r\in(0,\infty)}$ be the
same as in \eqref{teshu}.
Then there exists a function
$f\in C_{\rm{c}}^{\infty}(\mathbb{R})$ such that
$$\liminf_{r\to0^+}
\sup_{\lambda\in(0,\infty)}\lambda
\left[\nu_r(H_{\lambda,r}(f))\right]^\frac{1}{q}=\infty,$$
where $\nu_r$ and $H_{\lambda,r}$ are as, respectively,
in \eqref{nur} and \eqref{Hru}.
\end{theorem}

We now consider the BBM-type
formula involving the ball average operator
proposed in \cite[Remark 4.10]{jyy25}
[see also \eqref{equ-bbm-sim}].
The following theorem shows that
this formula also fails in general.
\begin{theorem}\label{pro-laplace}
Let $\theta\in(0,n)$,
$q\in(\max\{\frac{\theta}{2},1\},\infty)$,
and $\{\rho_r\}_{r\in(0,\infty)}$
be the same as in \eqref{teshu}.
Then, for any $f\in C_{\mathrm{c}}^3$,
\begin{equation}\label{equ-bbm}
\liminf_{r\to 0^+}
\int_{\mathbb{R}^n}
\int_{\mathbb{R}^n}
\frac{|f(y)-f_{B(y,r)}|^q}
{|x-y|^{2q}}\rho_{r}(x-y)\,dy\,dx<\infty
\end{equation}
if and only if $\Delta f\equiv 0$.
\end{theorem}
The next theorem gives the corresponding negative
result for the BSVY-type analogue involving ball averages.
\begin{theorem}\label{pro-weaknorm}
Let $\theta\in(0,n)$,
$q\in(\max\{\frac{\theta}{2},1\},\infty)$, and
$\{\rho_r\}_{r\in(0,\infty)}$ be the same as in \eqref{teshu}.
Then, for any $f\in C_{\mathrm{c}}^2$,
\begin{align}\label{equ-weakLq}
\liminf_{r\to 0^+}
\sup_{\lambda\in(0,\infty)}
\lambda
\left[\nu_{r}(E_{\lambda,r}(f))\right]^{\frac{1}{q}}<\infty
\end{align}
if and only if $f\equiv 0$, where $E_{\lambda,r}$ is defined in \eqref{Eru}.
\end{theorem}

In \cite[Remark 3.19(i)]{jyy25}, Jia et al. asked whether
a characterization of ball Banach Sobolev spaces involving
Taylor polynomials holds.
Theorem \ref{tailek} gives an affirmative answer to this question.
To state this result, we first recall
the definitions of the centered ball average
operators and the Hardy--Littlewood maximal operator.

For any $r\in(0,\infty)$, the
\textit{centered ball average operator} $\mathcal{B}_r$
is defined by setting, for any $f\in L^1_{\text{loc}}$
and $x\in\mathbb{R}^n$,
$$\mathcal{B}_r(f)(x):=\fint_{B(x,r)}f(y)\,dy.$$
We say that $\{\mathcal{B}_r\}_{r\in(0,\infty)}$
is \textit{uniformly bounded}
on the ball Banach function space $X$ if there
exists $C\in(0,\infty)$ such that,
for any $f\in X$,
$\sup_{r\in(0,\infty)}\|\mathcal{B}_r(f)\|_X\leq C\|f\|_X$.

Recall that the \textit{Hardy--Littlewood maximal operator}
$\mathcal{M}$ is defined by setting,
for any $f\in L^1_{\text{loc}}$ and $x\in\mathbb{R}^n$,
$$\mathcal{M}(f)(x):=\sup_{r\in(0,\infty)}
\fint_{B(x,r)}|f(y)|\,dy.$$
\begin{theorem}\label{tailek}
Let $k\in\mathbb{N}$, $q\in[1,\infty)$, and $X$ be a
ball Banach function space such that $X^{\frac{1}{q}}$
is also a ball Banach function space.
Suppose that $\mathcal{M}$ is bounded on $X$ and
$\{\mathcal{B}_r\}_{r\in(0,\infty)}$
are uniformly bounded on $X^{\frac{1}{q}}$. Further
assume that both $X$ and $X'$ have absolutely continuous
norms. Then $f\in W^{k,X}$ if and only if $f\in W^{k-1,X}$ and
there exists a function $H$ satisfying $H(x)\in\mathbb{M}_{n^{k-1},n}$
for almost every $x\in\mathbb{R}^n$ and
$|H(\cdot)|\in X$ such that
\begin{align}\label{equ-f-taylor}
\lim_{r\rightarrow0^+}\left\|\left[\fint_{B(\cdot,r)}
\frac{|f(y)-T_{\cdot}^{k-1}(f)(y)-\frac{1}{k!}
\sigma_k(y-\cdot)H(\cdot)(y-\cdot)|^q}{|y-\cdot|^{kq}}\,dy
\right]^{\frac{1}{q}}\right\|_X=0,
\end{align}
where $\mathbb{M}_{n^{k-1},n}$, $\sigma_k(\cdot)$,
and $|H(\cdot)|$ are as in Section \ref{s3} and $T_x^{k-1}(f)$ is
as in \eqref{taileduo}.
\end{theorem}
\begin{remark}
Under the assumptions of Theorem \ref{tailek},
Jia et al. proved that $f\in W^{k,X}$ if and only if $f\in W^{k-1,X}$
and the condition in \cite[(3.46)]{jyy25} holds.
In Theorem \ref{tailek}, we obtain a new characterization of $W^{k,X}$
by replacing the factor $\frac{1}{r^k}$ in \cite[(3.46)]{jyy25} with $\frac{1}{|y-x|^k}$
in \eqref{equ-f-taylor}, which complements \cite[Proposition 3.18]{jyy25}.
\end{remark}
In \cite[Remark 3.19(ii)]{jyy25}, Jia et al. asked whether
ball Banach Sobolev spaces can be characterized in terms of
higher-order minimizing polynomials [see \eqref{equ-local-van}].
The following theorem gives a negative answer to this question.
\begin{theorem}\label{the-4}
There exists $f\in L^2(\mathbb{R})\setminus
 W^{1,2}(\mathbb{R})$ satisfying
$$\lim_{r\to 0^+}
\left\|\frac{1}{r}
\left[\fint_{B(\cdot,r)}
\left|f(y)-P_{B(\cdot,r)}^{(1)}(f)(y)\right|^2\,dy\right]
^{\frac{1}{2}}
\right\|_{L^2(\mathbb{R})}=0.$$
\end{theorem}
Our final result concerns whether the
special John--Nirenberg--Campanato spaces via
congruent cubes $JN^{\rm{con}}_{(p,q,k)_\alpha}$ (see Definition \ref{def-john})
are independent of their
integrability exponent $q$ (see \cite[Remark 5.20]{jyy25} and
\cite[Remark 2.17]{jtyyz}).
\begin{theorem}\label{q-wug}
Let $p\in(1,\infty)$, $k\in\mathbb{N}$, and
$\alpha\in[0,\infty)$. Then, for any $q\in(1,p]$,
$JN^{\rm{con}}_{(p,q,k-1)_\alpha}=
JN^{\rm{con}}_{(p,1,k-1)_\alpha}$ with equivalent norms
if and only if $\alpha\neq0$.
\end{theorem}

In \cite[Remark 5.20]{jyy25}, Jia et al. asked whether
$JN_{(p,q,k-1)_{\alpha}}^{\mathrm{con}}=
JN_{(p,1,k-1)_{\alpha}}^{\mathrm{con}}$ when $k\in\mathbb{N}$,
$\alpha\in[0,\frac{k}{n}]$, $p\in(1,\infty)$,
and $q\in[1,p]$, whose special case $\alpha=\frac{k}{n}$
was proved to be true in \cite[Theorem 5.19]{jyy25}.
Theorem \ref{q-wug} gives a complete answer to this question
by covering the full range $\alpha\in[0,\infty)$.

Regarding the proofs of these theorems, the first key insight
is to combine orthogonal
polynomial projections with Plancherel's theorem to establish
an estimate for local oscillations in terms of the Fourier transform
(see Lemma \ref{zhengjiaotouying}). To construct the counterexample in
Theorem \ref{the-4}, a natural approach is to use
the first-order Taylor expansion to
estimate the remainder in terms of the second derivative.
However, this standard argument cannot be used here
because the desired counterexample requires
$\nabla f\notin L^2(\mathbb R)$, but $f\in L^2(\mathbb{R})$ and $\nabla^2f\in L^2(\mathbb{R})$
together imply $\nabla f\in L^2(\mathbb R)$, which is a contradiction.
To avoid this, we use the Fourier estimate
to control the normalized approximation error globally by choosing suitable
frequency components. This approach allows us to construct a function
in $L^2(\mathbb R)\setminus W^{1,2}(\mathbb{R})$
whose normalized approximation error tends to zero in $L^2(\mathbb{R})$.
The second key insight is to use a dyadic telescoping argument
(see Lemmas \ref{yl1}, \ref{yl2}, and \ref{yl3})
to establish the independence of the
integrability exponent $q$ for the special John--Nirenberg--Campanato spaces
$JN^{\rm{con}}_{(p,q,k)_{\alpha}}$ when $\alpha\not=0$. Recall that,
for the John--Nirenberg--Campanato
spaces, Tao et al. established a John--Nirenberg
inequality (see \cite[Theorem 4.3]{tyy19}), which was used
to prove that these spaces are independent of the integrability exponent $q$
(see \cite[Proposition 4.1]{tyy19}).
However, for $JN^{\rm{con}}_{(p,q,k)_{\alpha}}$, the
corresponding John--Nirenberg type inequality is still unknown.
The aforementioned approach allows us to obtain the desired independence without using
such an inequality. This, together with the counterexample for $\alpha=0$
[see \eqref{equ-JN-def-f}], proves Theorem \ref{q-wug}.

The remainder of this article is organized as follows.

Section \ref{s3} is devoted to investigating the sharpness of the
constants in the VMO-type characterization of ball Banach Sobolev spaces.
We prove that these constants are sharp uniformly over all admissible
ball Banach function spaces (see Theorem \ref{zjx})
and give an example showing that they need not be optimal
for a fixed space (see Proposition \ref{exp-3.4}).

In Section \ref{s4}, we disprove the BBM-type formula involving minimizing polynomials
and its BSVY-type analogue
(see, respectively, Theorems \ref{wenti1fouding} and \ref{ruofanshu1})
by constructing a family of singular radial kernels.

In Section \ref{s5}, we first establish a higher-order BBM-type
formula involving Taylor polynomials (Theorem \ref{taile}).
Applying this formula,
we establish a new VMO-type characterization of
ball Banach Sobolev spaces
(Theorem \ref{tailek}).

Section \ref{s7} is devoted to investigating the characterization
of ball Banach Sobolev spaces via $k$-th order minimizing polynomials.
We first establish a Fourier transform estimate for local oscillations
(Lemma \ref{zhengjiaotouying}) and construct
a counterexample showing that the corresponding vanishing oscillation
condition is insufficient to characterize
the Sobolev space (Theorem \ref{the-4}).

In Section \ref{s8}, we consider the BBM-type formula involving ball averages
and its BSVY-type analogue.
Using the  singular radial kernels
constructed in Section \ref{s4}, we prove that both of the formulas fail
(Theorems \ref{pro-laplace} and \ref{pro-weaknorm}).

Finally, in Section \ref{s9}, we investigate the dependence
of $JN^{\rm{con}}_{(p,q,k)_{\alpha}}$ on the integrability
exponent $q$. We prove that these spaces are independent of $q$ in the range
$q\in[1,p]$ when $\alpha\neq0$, and show that this independence fails
when $\alpha=0$ (see Theorem \ref{q-wug}).

We conclude this introduction with a few notational
conventions. Let $\mathbb{N}:=\{1,2,\ldots\}$,
$\mathbb{Z}_+:=\mathbb{N}\cup\{0\}$,
$\mathbb{Z}$
denote the set of all integers, and
$\mathbb{Z}^n:=\mathbb{Z}\times\cdots\times\mathbb{Z}$
(repeated $n$ times). We denote by the
\emph{symbol $C$} a positive constant which is
independent of the main parameters involved, but may vary
from line to line. The \emph{notation $f\lesssim g$}
means $f\leq Cg$. If $f\lesssim g$ and $g\lesssim f$,
we then write $f\sim g$. For any $s\in\mathbb{R}$,
we use the \emph{notation} $\lfloor s\rfloor$ to denote
the largest integer not greater than $s$.
All the functions
on $\mathbb{R}^n$ and subsets of $\mathbb{R}^n$
considered in this paper
are assumed to be Lebesgue measurable.
Given a set $E\subset\mathbb{R}^n$, we denote its
\emph{characteristic function} by $\mathbf{1}_E$,
the \emph{Lebesgue measure} by $|E|$,
and its \emph{cardinality} by $\# E$.
For any $\xi\in\mathbb{R}^n$ such that $|\xi|=1$,
$d\sigma(\xi)$ denotes the surface Lebesgue measure on the unit sphere
$\mathbb{S}^{n-1}\subset\mathbb{R}^n$ and $\sigma(\mathbb{S}^{n-1}):=\int_{\mathbb{S}^{n-1}}d\sigma$.
In addition, the \emph{space}
$\mathscr M$ denotes the set of all measurable
functions on ${\mathbb{R}^n}$.
For any
$p\in(0,\infty]$ and $E\subset\mathbb{R}^n$,
the \emph{Lebesgue space} $L^p(E)$ is defined
to be the set of all Lebesgue measurable
functions $f$ on $E$ such that
$\|f\|_{L^p(E)}:=[\int_E|f(x)|^p\,dx]^{\frac{1}{p}}<\infty$
with the usual modification made when $p=\infty$.
In this article, any \emph{cube} $Q$
has finite edge length
and all its edges parallel to the coordinate
axes, but $Q$ need not be open or closed;
moreover, we denote by the
\emph{notation $\ell(Q)$} the edge length of $Q$.
We also use $\mathbf{0}$ to denote the
origin of ${\mathbb{R}^n}$.
In addition, we need the  basic
inequality that, for any positive sequence
$\{\alpha_i\}_{i=1}^\infty$,
any $\theta\in[0,\infty)$, and any $N\in\mathbb{N}$,
\begin{equation}\label{basic-inequality}
\left(\sum_{i=1}^N{\alpha_i}\right)^{\theta}
\le\max\left\{1,N^{\theta -1}\right\}
\sum_{i=1}^N{\alpha_{i}^{\theta}},
\end{equation}
which, when $\theta\in[0,1]$, also holds for $N=\infty$.
Finally, in all proofs we consistently retain the notation
introduced in the original theorem (or related statement).

\section{Proof of Theorem \ref{zjx}}\label{s3}
In this section, we prove Theorem \ref{zjx}.
We begin by recalling the definition of
ball quasi-Banach function spaces
and ball Banach function spaces (see \cite[Definition 2.2]{shyy}).
\begin{definition}\label{def-ball-Banach}
A quasi-Banach function space $X\subset\mathscr{M}$,
equipped with a quasi-norm $\|\cdot\|_{X}$
which makes sense for all the
functions in $\mathscr{M}$, is called
a \emph{ball quasi-Banach function space}
if $X$ satisfies that
\begin{enumerate}
\item[(i)]for any $f\in\mathscr{M}$, if $||f||_X=0$,
then $f=0$ almost everywhere,
\item[(ii)]for any $f,g\in\mathscr{M}$, if $|g|\leq|f|$
almost everywhere, then $||g||_X\leq||f||_X$,
\item[(iii)]for any $\{f_m\}_{m\in\mathbb{N}}$
in $\mathscr{M}$
and $f\in\mathscr{M}$, if $0\leq f_m\uparrow f$
almost everywhere as $m\rightarrow\infty$, then
$||f_m||_X\uparrow ||f||_X$ as $m\rightarrow\infty$,
\item[(iv)]for any ball $B$ in $\mathbb{R}^n$,
$\mathbf{1}_{B}\in X$.
\end{enumerate}
Moreover, a ball
quasi-Banach function space $X$ is called
a \emph{ball Banach function space} if $X$
satisfies the following additional conditions:
\begin{enumerate}
\item[(v)]for any $f,g\in X$,
$||f+g||_X\leq||f||_X+||g||_X$,
\item[(vi)]for any ball $B$ in $\mathbb{R}^n$,
there exists a positive constant $C_{(B)}$,
depending on $B$,
such that, for any $f\in X$,
$\int_{B}|f(y)|\,dy\leq C_{(B)}||f||_X$.
\end{enumerate}
\end{definition}
\begin{remark}
\begin{enumerate}
\item[(i)]Let $X$ be a ball quasi-Banach function space
on $\mathbb{R}^n$. From \cite[Remark 2.6(i)]{yhyy}, we deduce that,
for any $f\in\mathscr{M}$,
$\|f\|_{X}=0$ if and only if $f=0$ almost everywhere.
\item[(ii)]By \cite[Theorem 2]{dfmn} and both (ii) and (iii) of Definition
\ref{def-ball-Banach}, we conclude that
every ball quasi-Banach function space is complete.
\item[(iii)]As mentioned in \cite[Remark 2.6(ii)]{yhyy}, we obtain an equivalent
formulation of Definition \ref{def-ball-Banach} by replacing
``any ball $B$'' with ``any bounded measurable set $E$''.
\item[(iv)]In Definition \ref{def-ball-Banach}, if we replace
``any ball $B$'' by ``any measurable set $E$ with finite measure'',
then we obtain the definition of (quasi-)Banach function spaces
introduced in \cite[Definitions 1.1 and 1.3]{bs}.
Hence, a (quasi-)Banach function space is always a ball
(quasi-)Banach function space and the converse is not necessarily true.
\item[(v)]In Definition 2.1, if we replace (iv) by the following
saturation property that, for any measurable set $E\subset\mathbb{R}^n$
with $|E|\in(0,\infty)$, there exists a measurable subset $F$
of $E$ with $|F|\in(0,\infty)$ such that $\mathbf{1}_F\in X$,
then we obtain the concept of quasi-Banach function spaces in Lorist and Nieraeth \cite{LN}.
Moreover, from \cite[Proposition 2.5]{zyy} (see also \cite[Proposition 4.21]{N}), it follows that, if
the quasi-normed vector space $X$ under consideration satisfies the additional assumption
that the Hardy-Littlewood maximal operator $\mathcal{M}$ is weakly bounded on $X$ or one of its
convexifications, then the definition of quasi-Banach function spaces in \cite{LN} coincides
with the definition of ball quasi-Banach function spaces. Thus, under this additional
assumption, working with ball quasi-Banach function spaces in the sense of Definition
\ref{def-ball-Banach} or quasi-Banach function spaces in the sense of \cite{LN} would yield precisely the same
results.
\end{enumerate}
\end{remark}
Next, we recall the definition of the
$p$-convexification for ball quasi-Banach function spaces,
which is exactly \cite[Definition 2.6(i)]{shyy}.
\begin{definition}\label{tuhua}
Let $X$ be a ball quasi-Banach function space and
$p\in(0,\infty)$.
The \textit{p-convexification} $X^p$ of $X$ is
defined by setting
$X^p:=\{f\in\mathscr{M}:|f|^p\in X\}$ equipped with the
\textit{quasi-norm}
$||f||_{X^p}:=|||f|^p||_{X}^{\frac{1}{p}}$
for any $f\in X^p$.
\end{definition}

In what follows, for any open set $U\subset\mathbb{R}^n$,
$\beta\in\mathbb{Z}_+^n$,
and $f\in L^1_{\text{loc}}(U)$, we say
$g\in L^1_{\text{loc}}(U)$ is the $\beta$-th
\textit{weak derivative} of
$f$ if, for any $\varphi\in C_{\text{c}}^{\infty}(U)$
(the set of all infinitely differentiable
functions on $U$
with compact support),
$$\int_{U}f(x)\partial^\beta\varphi(x)\,dx
=(-1)^{|\beta|}\int_{U}g(x)\varphi(x)\,dx,$$
here and thereafter,
for any $\beta:=(\beta_1,\ldots,\beta_n)\in\mathbb{Z}_+^n$,
$|\beta|:=\beta_1+\ldots+\beta_n$,
$\beta!:=\beta_1!\times\cdots\times\beta_n!$, and
$\partial^\beta:=(\frac{\partial}{\partial x_1})^{\beta_1}\cdots
(\frac{\partial}{\partial x_n})^{\beta_n}$.

We now introduce the definitions of ball Banach Sobolev spaces.

\begin{definition}\label{def-sobolev}
Let $k\in\mathbb{Z}_+$, $U$ be an open
subset of $\mathbb{R}^n$,
and $X$ be a ball Banach function space.
The \textit{homogeneous k-th order ball Banach Sobolev space}
$\dot{W}^{k,X}(U)$
is defined to be the set of
all  $f\in L^1_{\text{loc}}(U)$ satisfying
$|\nabla^kf|\in X(U)$,
equipped with the seminorm
$\|f\|_{\dot{W}^{k,X}(U)}:=\|\,|\nabla^kf|\,\|_{X(U)}
:=\|\,|\nabla^kf|\mathbf{1}_U\|_{X}$, where
\begin{equation}\label{def-wf}
\left|\nabla^kf\right|:=\left[\sum_{\{\beta\in\mathbb{Z}^n_+:
|\beta|=k\}}\frac{k!}{\beta!}
\left|\partial^\beta f(x)\right|^2\right]^{\frac{1}{2}}.	
\end{equation}
We denote $\dot{W}^{k,X}(\mathbb{R}^n)$ simply by $\dot{W}^{k,X}$.

The \textit{inhomogeneous k-th order ball Banach Sobolev space}
$W^{k,X}(U)$ is defined to be the set of
all  $f\in X(U)$ such that $\|f\|_{W^{k,X}(U)}
:=\sum_{d=0}^{k}\|\,\nabla^df\,\|_{X(U)}<\infty$.
We denote $W^{k,X}(\mathbb{R}^n)$ simply by $W^{k,X}$.

Moreover, for any $f\in L^1_{\text{loc}}$,
$f\in\dot{W}^{k,X}_{\text{loc}}$ means that,
for any bounded open subset
$U\subset\mathbb{R}^n$, $\|f\|_{\dot{W}^{k,X}(U)}<\infty$, and
$f\in W^{k,X}_{\text{loc}}$ means that,
for any bounded open subset
$U\subset\mathbb{R}^n$, $\|f\|_{W^{k,X}(U)}<\infty$.
\end{definition}
We now recall the definition of the minimizing polynomial.
For any $s\in\mathbb{Z}_+$, let $\mathcal{P}_s$ denote
the set of all polynomials of total degree not greater than $s$.
For any $k\in\mathbb{Z}_+$, any ball or any cube $B\subset \mathbb{R}^n$,
and any $f\in L^1_{\rm{loc}}$,
let $P^{(k)}_{B}$ denote the
\textit{k-th minimizing polynomial} of $f$ on $B$,
i.e., $P^{(k)}_{B}$ denotes the
unique polynomial $P\in\mathcal{P}_k$ such that
$\int_{B}[f(x)-P(x)]x^{\beta}\,dx=0$
for any $\beta \in\mathbb{Z}^n_+$ with $|\beta|\leq k$.
It is well known that, for any $f\in L^1_{\text{loc}}$
and any cube or any ball $B$ in $\mathbb{R}^n$, $P^{(0)}_{B}(f)=f_{B}:=\fint_{B}f(y)\,dy$,
and, for any $k\in\mathbb{N}$,
there exists a positive constant $C_{(k)}$
such that, for any $x\in B$,
\begin{equation}\label{jxdxsjfgj}
\left|P^{(k)}_{B}(f)(x)\right|\leq C_{(k)}\fint_{B}|f(y)|\,dy	
\end{equation}
(see, for instance, \cite[Lemma 4.11]{l}), which
holds with $B$ replaced by any cube $Q$.

We also recall the following elementary inequality
(see, for instance, \cite[(2.7)]{FMS(2018)}).
\begin{lemma}\label{fmsbds}
Let $q\in[1,\infty)$ and $\delta\in(0,1)$. Then, for any $a,b\in\mathbb{R}$,
$$|a+b|^q\leq (1+\delta)^q|a|^q+
\frac{(1+\delta)^q}{\delta^q}|b|^q.$$
\end{lemma}
To describe the constants in Theorem \ref{zjx}, we next
introduce some notation.
For any $n,m\in\mathbb{N}$
and any continuously
differentiable function $g:=(g_1,\ldots,g_m):\mathbb{R}^n
\rightarrow\mathbb{R}^m$, define
$$\nabla g:=
\left(\frac{\partial g_i}{\partial x_j}\right)
_{i\in\{1,\ldots,m\},j\in\{1,\ldots,n\}}
\in\mathbb{M}_{m,n},$$
where $\mathbb{M}_{m,n}$ denotes the set of all
$m\times n$ matrices.
For any $A=(a_{i,j})_{i\in\{1,\ldots,m\},j\in\{1,\ldots,n\}}
\in\mathbb{M}_{m,n}$,
let $|A|:=
(\sum_{i=1}^{m}\sum_{j=1}^{n}
|a_{i,j}|^2)^{\frac{1}{2}}$
and
$$\psi _{n,m}(A):=
\left(a_{11},\ldots,a_{1n},a_{21},\ldots,a_{2n},\ldots,
a_{m1},\ldots,a_{mn}\right).$$
Let $D^0(f):=f$, and,
for any $l\in\{1,\ldots,k\}$,
$D^l(f):=\psi _{n,n^{l-1}}(\nabla (D^{l-1}f))$.
For any $h\in\mathbb{R}^n$, let $\sigma_1(h):=1$ and,
for any
$l\in\mathbb{N}\cap[2,\infty)$,
$\sigma_l(h):=(h_1\sigma_{l-1}(h),
\ldots,h_n\sigma_{l-1}(h))$. By these, we have
\begin{align}\label{equ}
\left|\nabla\left(D^{k-1}f\right)(x)\right|
=\left|\nabla^kf(x)\right|.
\end{align}
\begin{remark}
We point out that there exists a gap in \cite[(3.4)]{jyy25}: with
\cite[Definition 2.6]{jyy25} which loses the correct
coefficients, the identity \eqref{equ},
which is precisely \cite[(3.4)]{jyy25}, does not hold in general.
For instance,
let $n=k=2$ and $f(x_1,x_2):=x_1x_2$. Then, for any $x=(x_1,x_2)\in\mathbb{R}^2$,
\begin{align*}
\left|\nabla\left(D^{k-1}f\right)(x)\right|^2&=
\left|\frac{\partial^2f}{\partial x_1^2}\right|^2+
\left|\frac{\partial^2f}{\partial x_1\partial x_2}\right|^2
+\left|\frac{\partial^2f}{\partial x_2\partial x_1}\right|^2
+\left|\frac{\partial^2f}{\partial x_2^2}\right|^2\\
&=2\neq1=\sum_{\{\beta\in\mathbb{Z}^n_+:
|\beta|=k\}}\left|\partial^\beta f(x)\right|^2.
\end{align*}
The definition in \eqref{def-wf} is the correct version of
\cite[Definition 2.6]{jyy25}, which
ensures that \eqref{equ} holds.
\end{remark}
To prove Theorem \ref{zjx}, we recall
the following matrix space introduced
in \cite[(3.5)]{jyy25}. For any $k\in\mathbb{N}$, let
\begin{align*}
\widetilde{\mathbb{M}}_{n^{k-1},n}
:=\left\{\nabla\left(D^{k-1}P\right):
P\in\mathcal{P}_k\setminus\mathcal{P}_{k-1}\right\}.
\end{align*}
The following lemma is precisely \cite[Theorem 3.3 and (3.12)]{jyy25}.
\begin{lemma}
Let $k\in\mathbb{N}$, $q\in[1,\infty)$, and $X$ be
a ball Banach function space such
that $X^{\frac{1}{q}}$ is also a ball Banach
function space. Then,
for any $f\in C^k$ with
$|\nabla^kf|\in C_{\rm{c}}$,
\begin{align}\label{equ-f-P}
\gamma{(n,q,k)}\left\|\,\left|\nabla^kf\right|\,
\right\|_X\leq
\lim_{r\rightarrow 0^+}\left\|\frac{1}{r^k}
\left[\fint_{B(\cdot,r)}
\left|f(y)-P^{(k-1)}_{B(\cdot,r)}(f)(y)
\right|^q\,dy\right]^{\frac{1}{q}}\right\|_X
\leq\Gamma{(n,q,k)}\left\|\,\left|\nabla^kf\right|
\,\right\|_X
\end{align}
and
\begin{align}\label{jixiaoduoxiangshijixian}
\lim_{r\rightarrow 0^+}
\left\|\frac{1}{r^k}\left[\fint_{B(\cdot,r)}
\left|f(y)-P^{(k-1)}_{B(\cdot,r)}(f)(y)\right|^q\,dy
\right]^{\frac{1}{q}}\right\|_X
=\left\|\mathrm{Osc}^{(q)}_k\left(\nabla\left(D^{k-1}f
\right)\right)\right\|_X,
\end{align}
where
\begin{equation}\label{gandGamma}
\gamma(n,q,k):=
\min_{\genfrac{}{}{0pt}{}{A\in
\widetilde{\mathbb{M}}_{n^{k-1},n}}
{|A|=1}}
\operatorname{Osc}_k^{(q)}(A),
\qquad
\Gamma(n,q,k):=
\max_{\genfrac{}{}{0pt}{}{A\in
\widetilde{\mathbb{M}}_{n^{k-1},n}}
{|A|=1}}
\operatorname{Osc}_k^{(q)}(A),
\end{equation}
and, for any $A\in \mathbb{M}_{n^{k-1},n}$,
\begin{equation}\label{OSCD}
\mathrm{Osc}^{(q)}_k(A):=\frac{1}{k!}
\left[\fint_{B(\mathbf{0},1)}
\left|\sigma_k(x)Ax-P^{(k-1)}_{B(\mathbf{0},1)}
(\sigma_k(\cdot)A\cdot)(x)\right|^q\,dx\right]^{\frac{1}{q}}.
\end{equation}
\end{lemma}

To prove Theorem \ref{zjx}, we establish the following result.
\begin{theorem}\label{wt1jd}
Let $k\in\mathbb{N}$ and $q\in[1,\infty)$.
Then there exist functions $f_{-},
f_{+}\in C_{\rm{c}}^{\infty}$ with the following property:
for any $\varepsilon\in(0,\infty)$,
there exist ball Banach function spaces
$X_\varepsilon$ and $Y_{\varepsilon}$,
with $X_\varepsilon^{\frac{1}{q}}$ and
$Y_\varepsilon^{\frac{1}{q}}$
also being
ball Banach function spaces, such that
\begin{align}\label{equ-gamma}
\displaystyle
\lim_{r\to0^+}
\left\|
\frac1{r^k}
\left[
\fint_{B(\cdot,r)}
\left|f_{-}(y)-P_{B(\cdot,r)}^{(k-1)}(f_{-})(y)\right|^q\,dy
\right]^{\frac{1}{q}}
\right\|_{X_\varepsilon}
\leq(\gamma(n,q,k)+\varepsilon)\left\|\ \left|\nabla^k f_{-}\right|\
\right\|_{X_\varepsilon}
\end{align}
and
\begin{align}\label{equ-Gamma}
\displaystyle
\lim_{r\to0^+}
\left\|\frac1{r^k}
\left[\fint_{B(\cdot,r)}
\left|f_{+}(y)-P_{B(\cdot,r)}^{(k-1)}(f_{+})(y)\right|^q\,dy
\right]^{\frac{1}{q}}
\right\|_{Y_\varepsilon}
\geq(\Gamma(n,q,k)-\varepsilon)\left\|\ \left|\nabla^k f_{+}\right|\ \right\|_{Y_\varepsilon}.
\end{align}
\end{theorem}
\begin{proof}
For any $r\in(0,\infty)$, $h\in L_{\mathrm{loc}}^1$,
and $x\in\mathbb{R}^n$, let
$$\mathcal T_r h(x)
:=
\frac1{r^k}
\left[
\fint_{B(x,r)}
\left|h(y)-P_{B(x,r)}^{(k-1)}(h)(y)\right|^q\,dy
\right]^{\frac{1}{q}}$$
and
$$\Theta_{n,k,q}:=\max_{\genfrac{}{}{0pt}{}{A\in
\mathbb{M}_{n^{k-1},n}}
{|A|=1}}
\operatorname{Osc}_k^{(q)}(A).$$
Since $\widetilde{\mathbb{M}}_{n^{k-1},n}\subset
\mathbb{M}_{n^{k-1},n}$
and $\mathbb{M}_{n^{k-1},n}$
is finite dimensional, it follows that
\begin{align}\label{jihejinji}
\Sigma:=\left\{A\in \widetilde{\mathbb{M}}_{n^{k-1},n}:
|A|=1\right\}
\end{align}
is compact.
This, together with the fact that
$\operatorname{Osc}_k^{(q)}$ is continuous
with respect to the matrix norm (see \cite[p.\,9]{jyy25}),
further implies that
there exist matrices $A_-,A_+\in\Sigma$ such that
$$\operatorname{Osc}_k^{(q)}(A_-)=\gamma(n,q,k)\ \mathrm{and}\
\operatorname{Osc}_k^{(q)}(A_+)=\Gamma(n,q,k).$$
By the definition of
$\widetilde{\mathbb{M}}_{n^{k-1},n}$,
we conclude that there exist polynomials
$Q_-,Q_+\in P_k\setminus P_{k-1}$ such that
\begin{align*}
\nabla\left(D^{k-1}Q_-\right)=A_-
\ \mathrm{and}\
\nabla\left(D^{k-1}Q_+\right)=A_+.
\end{align*}
By this, \eqref{equ}, and \eqref{jihejinji}, we obtain
\begin{align}\label{equ-nbla-Q}
\left|\nabla^kQ_-\right|
=\left|\nabla\left(D^{k-1}Q_-\right)\right|
=\left|A_-\right|=1=\left|A_+\right|
=\left|\nabla\left(D^{k-1}Q_+\right)\right|
=\left|\nabla^kQ_+\right|.
\end{align}
Choose $\eta\in C_{\rm{c}}^\infty$ satisfying
$\mathbf{1}_{B(\mathbf{0},2)}\le\eta\le
\mathbf{1}_{B(\mathbf{0},3)}$ and
define $f_-:=\eta Q_-$ and $f_+:=\eta Q_+$.
Then $f_-,f_+\in C_{\rm{c}}^\infty$.
From the fact that $f_{-}=Q_{-}$ and $f_{+}=Q_{+}$ on
$B(\mathbf{0},2)$,
we deduce that $\nabla^k{f_{-}}=\nabla^k{Q_{-}}$
and $\nabla^k{f_{+}}=\nabla^k{Q_{+}}$ on
$B(\mathbf{0},\frac{3}{2})$, which,
combined with \eqref{equ-nbla-Q} and
the fact that $Q_-,Q_+\in P_k\setminus P_{k-1}$,
further implies that
$$\left\| \ \left|\nabla ^kf_{-}\right|\
\right\|_{L^q[\mathbb{R}^n\setminus B(\mathbf{0},1)]}
\geq \left\| \ \left|\nabla^kQ_{-}\right| \
\right\|_{L^q[B(\mathbf{0},\frac{3}{2})\setminus
B(\mathbf{0},1)]}
=\left|B\left(\mathbf{0},\frac{3}{2}\right)\setminus
B(\mathbf{0},1)
\right|^{\frac{1}{q}}>0$$
and
$$\left\| \ \left|\nabla ^kf_{+}\right|\
\right\|_{L^q[\mathbb{R}^n\setminus B(\mathbf{0},1)]}\geq
\left\| \ \left|\nabla ^kQ_{+}\right|\
\right\|_{L^q[B(\mathbf{0},\frac{3}{2})\setminus
B(\mathbf{0},1)]}
=\left|B\left(\mathbf{0},\frac{3}{2}\right)\setminus
B(\mathbf{0},1)
\right|^{\frac{1}{q}}>0.$$
On the other hand, since
$f_-,f_+\in C_{\mathrm{c}}^{\infty}$,
it follows that
$|\nabla^kf_-|,|\nabla^kf_+|\in L^q$.
Fix $\varepsilon\in(0,\infty)$, define
$$w_{\varepsilon}^{X}:=\mathbf{1}_{B(\mathbf{0},1)}+
\frac{|B(\mathbf{0},1)|{\varepsilon}^{q}}
{\Theta_{n,k,q}^q
|| \ |\nabla ^kf_{-}|\ ||_{L^q[\mathbb{R}^n\setminus
B(\mathbf{0},1)]}^q}
\mathbf{1}_{\mathbb{R}^n\setminus B(\mathbf{0},1)},$$
and let
$X_\varepsilon:=L^q_{w_\varepsilon^{X}}$
be the set of all $h\in L_{\rm{loc}}^q$ satisfying
$$\|h\|_{X_\varepsilon}:=\left[
\int_{\mathbb R^n}|h(x)|^q w_\varepsilon^{X}(x)\,dx
\right]^{\frac{1}{q}}<\infty.$$
Similarly, let
$$w_{\varepsilon}^{Y}:=\mathbf{1}_{B(\mathbf{0},1)}+
\frac{|B(\mathbf{0},1)
|{\varepsilon}^{q}}{[\Gamma(n,q,k)]^{q}||
\ |\nabla^kf_{+}| \ ||_{L^q[\mathbb{R}^n\setminus
B(\mathbf{0},1)]}^q}
\mathbf{1}_{\mathbb{R}^n\setminus B(\mathbf{0},1)},$$
and let
$Y_\varepsilon:=L^q_{w_\varepsilon^{Y}}$
be the set of all $h\in L_{\rm{loc}}^q$ satisfying
$$\|h\|_{Y_\varepsilon}:=\left[
\int_{\mathbb R^n}|h(x)|^q w_\varepsilon^{Y}(x)\,dx
\right]^{\frac{1}{q}}<\infty.$$
Then it is easy to verify that
$X_\varepsilon, Y_{\varepsilon},
X_\varepsilon^{\frac{1}{q}},$ and
$Y_{\varepsilon}^{\frac{1}{q}}$
are all ball Banach function spaces.
By \eqref{jixiaoduoxiangshijixian}, we obtain
\begin{align*}
\lim_{r\rightarrow 0^+}
\frac{||\mathcal{T}_rf_{-}||_{X_{\varepsilon}}}
{||\ |\nabla^kf_{-}|\ ||_{X_{\varepsilon}}}
=\frac{||\text{Osc}_{k}^{(q)}
(\nabla(D^{k-1}f_{-}))||_{X_{\varepsilon}}}
{||\ |\nabla^kf_{-}|\ ||_{X_{\varepsilon}}}.
\end{align*}
On the other hand, from \eqref{basic-inequality} and
the definition of $X_{\varepsilon}$,
we infer that
\begin{align*}
&\left\|\text{Osc}_{k}^{(q)}
\left(\nabla\left(D^{k-1}f_{-}\right)\right)
\right\|_{X_{\varepsilon}}\\
&\quad
\leq\left\|\text{Osc}_{k}^{(q)}
\left(\nabla\left(D^{k-1}f_{-}\right)\right)
\right\|_{L^q[B(\mathbf{0},1)]}+
\frac{|B(\mathbf{0},1)|^{\frac{1}{q}}
\varepsilon||\text{Osc}_{k}^{(q)}(\nabla(D^{k-1}f_{-}))
||_{L^q[\mathbb{R}^n\setminus
B(\mathbf{0},1)]}}{\Theta_{n,k,q}
||\ |\nabla ^kf_{-}|\ ||_{L^q[\mathbb{R}^n\setminus
B(\mathbf{0},1)]}}\\
&\quad
\leq
\left\|\text{Osc}_{k}^{(q)}
\left(\nabla\left(D^{k-1}Q_{-}\right)\right)
\right\|_{L^q[B(\mathbf{0},1)]}
+\frac{|B(\mathbf{0},1)|^{\frac{1}{q}}
\varepsilon\Theta_{n,k,q}
||\ |\nabla ^kf_{-}|\ ||_{L^q[\mathbb{R}^n
\setminus B(\mathbf{0},1)]}}
{\Theta_{n,k,q}||\ |\nabla ^kf_{-}|\
||_{L^q[\mathbb{R}^n\setminus B(\mathbf{0},1)]}}\\
&\quad=\gamma(n,q,k)|B(\mathbf{0},1)|^{\frac{1}{q}}+
\varepsilon|B(\mathbf{0},1)|^{\frac{1}{q}}.
\end{align*}
This, together with the fact that
$|\nabla^kf_-|=|\nabla^kQ_-|\equiv1$ on $B(\mathbf{0},1)$,
further implies that
\begin{align*}
\lim_{r\rightarrow 0^+}
\frac{||\mathcal{T}_rf_{-}||_{X_{\varepsilon}}}
{||\ |\nabla^kf_{-}|\ ||_{X_{\varepsilon}}}
\leq
\frac{\gamma(n,q,k)|B(\mathbf{0},1)|^{\frac{1}{q}}
+\varepsilon|B(\mathbf{0},1)|^{\frac{1}{q}}}
{||\ |\nabla^kf_{-}|\ ||_{L^q[B(\mathbf{0},1)]}}
=\gamma(n,q,k)+\varepsilon,
\end{align*}
which shows \eqref{equ-gamma}. Next, we prove
\eqref{equ-Gamma}.
Indeed, by \eqref{jixiaoduoxiangshijixian}
and the definition of $Y_{\varepsilon}$, we conclude that
\begin{align}\label{equ-f-plus}
\lim_{r\rightarrow 0^+}
\frac{||\mathcal{T}_rf_{+}||_{Y_{\varepsilon}}}
{||\ |\nabla^kf_{+}|\ ||_{Y_{\varepsilon}}}
=\frac{||\text{Osc}_{k}^{(q)}
(\nabla(D^{k-1}f_{+}))||_{Y_{\varepsilon}}}
{||\ |\nabla^kf_{+}|\ ||_{Y_{\varepsilon}}}
\geq \frac{||\text{Osc}_{k}^{(q)}(\nabla(D^{k-1}f_{+}))
||_{L^q[B(\mathbf{0},1)]}}
{||\ |\nabla^kf_{+}|\ ||_{Y_{\varepsilon}}}.
\end{align}
From the definition of $Y_{\varepsilon}$,
we deduce that
\begin{align*}
\left\|\ \left|\nabla^kf_{+}\right|\
\right\|_{Y_{\varepsilon}}
&\leq\left\|\ \left|\nabla^kf_{+}\right|\
\right\|_{L^q[B(\mathbf{0},1)]}
+\frac{|B(\mathbf{0},1)|^{\frac{1}{q}}\varepsilon}
{\Gamma(n,q,k)
||\ |\nabla ^kf_{+}|\ ||_{L^q[\mathbb{R}^n\setminus
B(\mathbf{0},1)]}}
\left\|\ \left|\nabla^kf_{+}\right|\
\right\|_{L^q[\mathbb{R}^n\setminus B(\mathbf{0},1)]}\\
&=|B(\mathbf{0},1)|^{\frac{1}{q}}+
\frac{|B(\mathbf{0},1)|^{\frac{1}{q}}\varepsilon}
{\Gamma(n,q,k)}.
\end{align*}
This, combined with \eqref{equ-f-plus}, further implies that
\begin{align*}
\lim_{r\rightarrow 0^+}
\frac{||\mathcal{T}_rf_{+}||_{Y_{\varepsilon}}}
{||\ |\nabla^kf_{+}|\ ||_{Y_{\varepsilon}}}
\geq\frac{\Gamma(n,q,k)
|B(\mathbf{0},1)|^{\frac{1}{q}}}
{|B(\mathbf{0},1)|^{\frac{1}{q}}+
\frac{|B(\mathbf{0},1)|^{\frac{1}{q}}\varepsilon}
{\Gamma(n,q,k)}}
=\frac{\Gamma(n,q,k)}
{1+\frac{\varepsilon}{\Gamma(n,q,k)}}
\geq\Gamma(n,q,k)-\varepsilon,
\end{align*}
which completes the proof of Theorem \ref{wt1jd}.
\end{proof}
We now apply Theorem \ref{wt1jd} to prove Theorem \ref{zjx}.
\begin{proof}[Proof of Theorem \ref{zjx}]
Suppose that there exists $\gamma'(n,q,k)\in
(\gamma(n,q,k),\infty)$
such that \eqref{equ-f-P} holds with $\gamma(n,q,k)$
replaced by $\gamma'(n,q,k)$. Choose $\varepsilon
\in(0,\frac{\gamma'(n,q,k)-\gamma(n,q,k)}{2})$,
and let $X_{\varepsilon}$ and $f_-$ be the same as in
\eqref{equ-gamma}. By these, \eqref{equ-f-P},
and \eqref{equ-gamma}, we obtain
\begin{align*}
\gamma'(n,q,k)&\le
\lim_{r\to0^+}
\frac{\|\frac1{r^k}[\fint_{B(\cdot,r)}
|f_-(y)-P_{B(\cdot,r)}^{(k-1)}(f_-)(y)|^q\,dy
]^{\frac{1}{q}}
\|_{X_\varepsilon}}{\|\ |\nabla^k f_{-}|\ \|_{X_\varepsilon}}\\
&\le\gamma(n,q,k)+\frac{\gamma'(n,q,k)-\gamma(n,q,k)}{2}
=\frac{\gamma'(n,q,k)+\gamma(n,q,k)}{2},
\end{align*}
which contradicts the assumption that
$\gamma'(n,q,k)>\gamma(n,q,k)$.
Applying an argument similar to the above
with $\gamma(n,q,k)$, $f_-$, and $X_{\varepsilon}$,
respectively, replaced by $\Gamma(n,q,k)$, $f_+$,
and $Y_{\varepsilon}$,
we also show that $\Gamma(n,q,k)$ is sharp.
This completes the proof of Theorem \ref{zjx}.
\end{proof}
Theorem \ref{zjx} shows that $\gamma(n,q,k)$ and $\Gamma(n,q,k)$
are sharp uniformly over all the ball Banach function spaces
satisfying the hypotheses therein.
However, they need not be sharp for a fixed space $X$.
The following proposition shows this for
$X=L^2(\mathbb{R}^2)$.
\begin{proposition}\label{exp-3.4}
Let $f\in C^2(\mathbb{R}^2)$ and $|\nabla^2f|\in C_{\rm{c}}(\mathbb{R}^2)$. Then
$$\lim_{r\rightarrow0^+}\left\|\frac{1}{r^2}
\left[\fint_{B(\cdot,r)}
\left|f(y)-P^{(1)}_{B(\cdot,r)}(f)(y)
\right|^2\,dy\right]^{\frac{1}{2}}\right\|_{L^2(\mathbb{R}^2)}
=\frac{1}{8}\left\|\,\left|\nabla^2f\right|
\,\right\|_{L^2(\mathbb{R}^2)}.$$
\end{proposition}
To prove Proposition \ref{exp-3.4}, we
recall more concepts.
Let $\mathcal{S}$ denote the set of all \emph{Schwartz functions} on $\mathbb{R}^n$
and $\mathcal{S}'/\mathcal{P}$ the set of all tempered
distributions modulo polynomials. Recall that,
for any $f\in L^1$, its \emph{Fourier transform} $\mathcal{F}f$
is defined by setting, for any $\xi\in\mathbb{R}^n$,
\begin{align*}
\mathcal{F}f(\xi):=
(2\pi)^{-\frac{n}{2}}
\int_{\mathbb{R}^n}f(x)e^{-ix\cdot\xi}\,dx
\end{align*}
and its \emph{inverse Fourier transform} is defined by setting
$\mathcal{F}^{-1}f(\xi):=\mathcal{F}f(-\xi)$ for any $\xi\in\mathbb{R}^n$.
These definitions are naturally extended to $\mathcal{S}'/\mathcal{P}$.

We now prove Proposition \ref{exp-3.4}.
\begin{proof}[Proof of Proposition \ref{exp-3.4}]
For any $A=(a_{i,j})_{i\in\{1,2\},j\in\{1,2\}}\in\widetilde{\mathbb{M}}_{2,2}
\subset\{A\in\mathbb{M}_{2,2}:A=A^{T}\}$, let
$\text{tr}(A):=a_{11}+a_{22}$. Then we have
$P_{B(\mathbf{0},1)}^{(1)}(\sigma_2(\cdot)A\cdot)=\frac{\rm{tr}A}{4}.$
This implies that
\begin{align}\label{oscgj}
\mathrm{Osc}^{(2)}_2(A)&=\frac{1}{2}
\left[\fint_{B(\mathbf{0},1)}
\left|x^{T}Ax-\frac{\rm{tr}A}{4}
\right|^2\,dx\right]^{\frac{1}{2}}=\frac{1}{4\sqrt{3}}\left(
|A|^2-\frac{|\rm{tr}A|^2}{4}\right)^{\frac{1}{2}}.
\end{align}
For any $f\in C^2$ with $|\nabla^2f|\in C_{\rm{c}}$,
we conclude that, for any $\beta\in\mathbb{Z}_+^2$ with $|\beta|=2$,
$$\mathcal{F}\left(\partial^\beta f\right)=-\xi^\beta\mathcal{F}f
=\frac{\xi^\beta}{|\xi|^2}\mathcal{F}(\Delta f).$$
By this, the fact that $|\nabla^2f|\in C_{\rm{c}}\subset L^1\cap L^2$,
and Plancherel's theorem, we obtain
$$\left\|\,\left|\nabla^2f\right|\,\right\|_{L^2(\mathbb{R}^2)}^2=
\|\Delta f\|_{L^2(\mathbb{R}^2)}^2.$$
From this and \eqref{oscgj}, we deduce that
\begin{align*}
\left\|\mathrm{Osc}^{(2)}_2\left(\nabla\left(D^1f
\right)\right)\right\|_{L^2(\mathbb{R}^2)}^2&
=\frac{1}{48}\left\|
\left|\nabla^2f\right|^2-\frac{1}{4}\left|
\frac{\partial^2f}{\partial x_1^2}+
\frac{\partial^2f}{\partial x_2^2}\right|^2
\right\|_{L^1(\mathbb{R}^2)}\\
&=\frac{1}{48}\left[\left\|\,\left|\nabla^2f\right|\,\right\|_{L^2(\mathbb{R}^2)}^2
-\frac{1}{4}\|\Delta f\|_{L^2(\mathbb{R}^2)}^2\right]
=\frac{1}{64}\left\|\,\left|\nabla^2f\right|\,\right\|_{L^2(\mathbb{R}^2)}^2.
\end{align*}
This completes the proof of Proposition \ref{exp-3.4}.
\end{proof}
\begin{remark}
\begin{enumerate}
\item[(i)]We claim that, for a fixed ball Banach function space $X$,
the positive constants $\gamma(n,q,k)$ and $\Gamma(n,q,k)$ in Theorem \ref{zjx}
need not be sharp.
For example, let $n=q=k=2$ and $X=L^2(\mathbb{R}^2)$.
Proposition \ref{exp-3.4} shows that the sharp lower and upper
constants are both $\frac{1}{8}$. To compute $\gamma(2,2,2)$ and $\Gamma(2,2,2)$,
let
$$|\rm{tr}A_1|^2:=\left|\rm{tr}\left(\frac{1}{\sqrt{2}}\rm{diag}
\{1,-1\}\right)\right|^2=0\leq|\rm{tr}A|^2\leq2|A|^2=2
=\left|\rm{tr}\left(\frac{1}{\sqrt{2}}\rm{diag}
\{1,1\}\right)\right|^2=:|\rm{tr}A_2|^2.$$
By \eqref{oscgj}, we obtain
\begin{equation*}
\gamma(2,2,2)=\min_{\genfrac{}{}{0pt}{}{A\in
\widetilde{\mathbb{M}}_{2,2}}{|A|=1}}
\operatorname{Osc}_2^{(2)}(A)=\min_{\genfrac{}{}{0pt}{}{A\in
\widetilde{\mathbb{M}}_{2,2}}{|A|=1}}\frac{1}{4\sqrt{3}}\left(
1-\frac{|\rm{tr}A|^2}{4}\right)^{\frac{1}{2}}=\operatorname{Osc}_2^{(2)}(A_2)
=\frac{1}{4\sqrt{6}}
\end{equation*}
and
\begin{equation*}
\Gamma(2,2,2)=\max_{\genfrac{}{}{0pt}{}{A\in
\widetilde{\mathbb{M}}_{2,2}}{|A|=1}}
\operatorname{Osc}_2^{(2)}(A)=\max_{\genfrac{}{}{0pt}{}{A\in
\widetilde{\mathbb{M}}_{2,2}}{|A|=1}}\frac{1}{4\sqrt{3}}\left(
1-\frac{|\rm{tr}A|^2}{4}\right)^{\frac{1}{2}}=\operatorname{Osc}_2^{(2)}(A_1)
=\frac{1}{4\sqrt{3}}.
\end{equation*}
Thus, we have
$$\gamma(2,2,2)<\frac{1}{8}<\Gamma(2,2,2),$$
which proves the claim.
\item[(ii)]For a fixed ball Banach function space $X$,
we define the following sharp constants:
$$\gamma(n,q,k;X):=\inf_{\genfrac{}{}{0pt}{}{f\in C^k,|\nabla^kf|\in C_{\rm{c}}}
{\|\,|\nabla^kf|\,\|_{X}\neq0}}
\frac{\|\mathrm{Osc}^{(q)}_k(\nabla(D^{k-1}f))\|_X}
{\|\,|\nabla^kf|\,\|_X}$$
and
$$\Gamma(n,q,k;X):=\sup_{\genfrac{}{}{0pt}{}{f\in C^k,|\nabla^kf|\in C_{\rm{c}}}
{\|\,|\nabla^kf|\,\|_{X}\neq0}}
\frac{\|\mathrm{Osc}^{(q)}_k(\nabla(D^{k-1}f))\|_X}
{\|\,|\nabla^kf|\,\|_X}.$$
It is interesting to know how to characterize the sharp constants
$\gamma(n,q,k;X)$ and $\Gamma(n,q,k;X)$ for a given
ball Banach function space $X$.
\end{enumerate}
\end{remark}

\section{Proofs of Theorems \ref{wenti1fouding} and
\ref{ruofanshu1}}\label{s4}
This section is divided into two subsections.
In Subsection \ref{sub-min-BBM}, we disprove the
proposed BBM-type formula involving minimizing polynomials
(see Theorem \ref{wenti1fouding}). Later, we show that its
BSVY-type analogue also fails (see Theorem \ref{ruofanshu1})
in Subsection \ref{sub-weak-min}.

\subsection{Proof of Theorem \ref{wenti1fouding}}\label{sub-min-BBM}
In this subsection, we prove Theorem \ref{wenti1fouding}.
To this end, we first recall the question in
\cite[Remark 3.14]{jyy25}. Let $X$
be a ball Banach function space satisfying
some mild assumptions (see \cite[Theorem 3.3]{jyy25}).
The question is whether, for any $q\in[1,\infty)$, $k\in\mathbb{N}$,
and $f\in C^k$ with $|\nabla^k f|\in C_{\mathrm{c}}$,
\begin{align*}
\lim_{r\to 0^+}
\left\|\left[\int_{\mathbb{R}^n}
\frac{|f(y)-P_{B(\cdot,r)}^{(k-1)}(f)(y)|^q}{|\cdot-y|^{kq}}
\rho_r(\cdot-y)\,dy\right]^{\frac{1}{q}}\right\|_{X}
\sim\left\|\ \left|\nabla^kf\right|\ \right\|_{X},
\end{align*}
where $\{\rho_r\}_{r\in(0,\infty)}$ satisfies the conditions
of Definition \ref{djh} and where
the positive equivalence constants are
independent of $f$.
\begin{definition}\label{djh}
A family $\{\rho_r\}_{r\in(0,\infty)}$ of nonnegative locally integrable functions
is called an \emph{approximation of the identity of radial type on} $\mathbb{R}^n$
if the following conditions hold:	
\begin{enumerate}
\item [$\rm{(i)}$] $\rho_r$ is a radial
function for any $r\in(0,\infty)$,
\item [$\rm{(ii)}$] $\int_{\mathbb{R}^n}\rho_{r}(x)\,
dx=1$ for any $r\in(0,\infty)$,
\item [$\rm{(iii)}$] For any $t\in(0,\infty)$,
$\lim_{r\to 0^+}\int_{\mathbb{R}^n\setminus B(\mathbf{0},t)}
\rho_r(x)\,dx=0$.
\end{enumerate}	
\end{definition}
Theorem \ref{wenti1fouding} gives a negative answer to the aforementioned question.
Specifically, let $q\in[1,\infty)$, $k:=1$, and
$X:=L^q$. Then, by the definition of the minimizing polynomial,
for any $f\in L_{\mathrm{loc}}^1$
and any ball $B$ in $\mathbb{R}^n$, $P_{B}^{(0)}(f)=f_B$.

To prove Theorem \ref{wenti1fouding},
we need the following lemma, which
is essentially \cite[Theorem 2.7]{dgpyyz24}.
\begin{lemma}\label{dgpyyz}
Let $q\in[1,\infty)$,
$\{\rho_{\varepsilon}\}_{\varepsilon\in(0,\infty)}$
satisfy the conditions of Definition \ref{djh}, and $X$ be a
ball Banach function space.
Assume that,
for any given $M\in(0,\infty)$, there exists
$N\in(0,\infty)$
such that
\begin{align*}
\lim_{\varepsilon \to 0+}
\left\|\left[\int_{B(\mathbf{0},M)}
\frac{\rho_{\varepsilon }(\cdot-y)}
{|\cdot-y|^q}\,dy\right]^{\frac{1}{q}}
\mathbf{1}_{\mathbb{R}^n\setminus
B(\mathbf{0},N)}\right\|_{X}=0.
\end{align*}
Then there exists a positive constant $C_{(q,n)}$ such that,
for any $f\in C_{\rm{c}}^2$,
\begin{align*}
\lim_{\varepsilon\to 0^+}
\left\|\left[\int_{\mathbb{R}^n}
\frac{|f(y)-f(\cdot)|^q}{|\cdot-y|^{q}}
\rho_\varepsilon(\cdot-y)\,dy\right]^{\frac{1}{q}}
\right\|_{X}
=C_{(q,n)}\left\|\ \left|\nabla f\right|\ \right\|_{X}.
\end{align*}
\end{lemma}
We now prove Theorem \ref{wenti1fouding}.
\begin{proof}[Proof of Theorem \ref{wenti1fouding}]
First, we verify that $\{\rho_r\}_{r\in(0,\infty)}$
satisfies the conditions of Definition \ref{djh}.
Condition (i) is evident, and hence we only need to
show (ii) and (iii). Fix $r\in(0,\infty)$.
Using the assumption that $\theta>0$, we conclude that
\begin{align}\label{equ-rho-1}
\int_{\mathbb{R}^n}
\rho_r(x)\,dx
&=\frac{\theta}{\sigma(\mathbb{S}^{n-1})r^\theta}
\int_{B(\mathbf{0},r)}|x|^{\theta-n}\,dx=
\frac{\theta}{r^\theta}
\int_{0}^{r}\lambda^{\theta-1}\,d\lambda=1,
\end{align}
which shows (ii). For any $\delta\in(0,\infty)$,
if $r\in(0,\delta )$, then $\mathop{\mathrm{supp}}(\rho_r)\subset
\overline{B(\mathbf{0},r)}\subset B(\mathbf{0},\delta)$, which
further implies that
\begin{equation}\label{equ-rho-0}
\int_{\mathbb{R}^n\setminus B(\mathbf{0},\delta)}
\rho_{r}(x)\,dx=0,
\end{equation}
and hence
$\lim_{r\to0^+}\int_{\mathbb{R}^n\setminus
B(\mathbf{0},\delta)}{\rho_{r}(x)}\,dx=0$,
which verifies (iii).

Next, we show Theorem \ref{wenti1fouding}.
The sufficiency is obvious and hence we only need
to show the necessity.
Fix $f\in C_{\mathrm{c}}^{2}$ satisfying \eqref{equ-infty},
and assume $f\not\equiv0$.
For any $M\in(0,\infty)$, $N\in(M+1,\infty)$, and
$r\in(0,1)$, if
$y\in B(\mathbf{0},M)$ and $x\in \mathbb{R}^n\setminus
B(\mathbf{0},N)$, then
$|x-y|\geq 1>r$ and hence $\rho_{r}(x-y)=0$.
This further implies that
$$\lim_{r \to 0+}
\left\|\left[\int_{B(\mathbf{0},M)}
\frac{\rho_{r}(\cdot-y)}
{|\cdot-y|^q}\,dy\right]^{\frac{1}{q}}
\mathbf{1}_{\mathbb{R}^n\setminus B(\mathbf{0},N)}
\right\|_{L^q}=0.$$
By this and Lemma \ref{dgpyyz} with $X:=L^q$, we conclude that
\begin{align}\label{yzjx}
\lim_{r\to 0^+}
\left\|
\left[\int_{\mathbb{R}^n}
\frac{|f(y)-f(\cdot)|^q}{|\cdot-y|^{q}}
\rho_r(\cdot-y)\,dy\right]^{\frac{1}{q}}
\right\|_{L^q}
=C_{(q,n)}\left\|\ \left|\nabla f\right|\ \right\|_{L^q}.
\end{align}	
On the other hand, we have
\begin{align*}
&\left[\int_{\mathbb{R}^n}
\int_{\mathbb{R}^n}
\frac{|f(x)-f_{B(x,r)}|^q}{|x-y|^q}
\rho_r(|x-y|)\,dy\,dx\right]^{\frac{1}{q}}\\
&\quad\lesssim \left[\int_{\mathbb{R}^n}
\int_{\mathbb{R}^n}
\frac{|f(y)-f_{B(x,r)}|^q}{|x-y|^q}
\rho_r(|x-y|)\,dy\,dx\right]^{\frac{1}{q}}+
\left[\int_{\mathbb{R}^n}
\int_{\mathbb{R}^n}
\frac{|f(x)-f(y)|^q}{|x-y|^q}
\rho_r(|x-y|)\,dy\,dx\right]^{\frac{1}{q}}.
\end{align*}
From \eqref{yzjx}, we deduce that
\begin{align}\label{equ-liminf}
&\liminf_{r\to 0^+}\left[\int_{\mathbb{R}^n}
\int_{\mathbb{R}^n}
\frac{|f(x)-f_{B(x,r)}|^q}{|x-y|^q}
\rho_r(x-y)\,dy\,dx\right]^{\frac{1}{q}}\notag\\
&\quad\lesssim \liminf_{r\to 0^+}\left[\int_{\mathbb{R}^n}
\int_{\mathbb{R}^n}
\frac{|f(y)-f_{B(x,r)}|^q}{|x-y|^q}
\rho_r(x-y)\,dy\,dx\right]^{\frac{1}{q}}+\left\|\ \left|
\nabla f\right|\ \right\|_{L^q(\mathbb{R}^n)}.
\end{align}
We claim that
\begin{align}\label{equ-f-B}
\liminf_{r\to 0^+}\left[\int_{\mathbb{R}^n}
\int_{\mathbb{R}^n}
\frac{|f(x)-f_{B(x,r)}|^q}{|x-y|^q}
\rho_r(x-y)\,dy\,dx\right]^{\frac{1}{q}}=\infty.
\end{align}
If this claim holds, then, from
\eqref{equ-liminf} and the fact that
$|\nabla f|\in L^q$, we infer that
\begin{align*}
\liminf_{r\to 0^+}\left[\int_{\mathbb{R}^n}
\int_{\mathbb{R}^n}
\frac{|f(y)-f_{B(x,r)}|^q}{|x-y|^q}
\rho_r(x-y)\,dy\,dx\right]^{\frac{1}{q}}=\infty,
\end{align*}
which contradicts \eqref{equ-infty} and further implies that
$f\equiv0$.
Therefore, it suffices to show \eqref{equ-f-B}.
Fix $r\in(0,1)$.
Indeed, by a change of variables, we obtain
\begin{align}\label{equ-rho-*}
\left[\int_{\mathbb{R}^n}
\int_{\mathbb{R}^n}
\frac{|f(x)-f_{B(x,r)}|^q}{|x-y|^q}
\rho_r(x-y)\,dy\,dx\right]^{\frac{1}{q}}=
\left[\int_{\mathbb{R}^n}
|f(x)-f_{B(x,r)}|^q\int_{\mathbb{R}^n}
\frac{\rho_r(y)}{|y|^q}
\,dy\,dx\right]^{\frac{1}{q}}.
\end{align}
From the definition of $\{\rho_r\}_{r\in(0,\infty)}$
and the assumption $\theta<q$, we deduce that
\begin{align*}
\int_{\mathbb{R}^n}
\frac{\rho_r(y)}{|y|^q}
\,dy\sim\int_{B(\mathbf{0},r)}
\frac{\theta }{r^\theta |y|^{n+q-\theta }}
\,dy=\infty.
\end{align*}
By this and \eqref{equ-rho-*},
we only need to show that
\begin{align}\label{equ-f-fav}
\int_{\mathbb{R}^n}
|f(x)-f_{B(x,r)}|^q\,dx>0.
\end{align}
If this does not hold, then,
by the continuity of $f$ and $f_{B(\cdot,r)}$, we
conclude that
$f\equiv f_{B(\cdot,r)}$.
Choose $x_0\in\mathbb{R}^n$ such that
$f(x_0)=\min_{x\in\mathbb{R}^n}f(x)=:m$.
This, together with $f(x_0)=f_{B(x_0,r)}$,
further implies that
$f\equiv m$ on $B(x_0,r)$. Therefore,
for any $y\in B(x_0,r)$, we have $f(y)=f_{B(y,r)}=m$,
which further implies that $f\equiv m$ on $B(y,r)$. On the other hand,
a geometric observation shows that
\begin{align*}
\bigcup_{y\in B(x_0,r)}B(y,r)=B(x_0,2r).
\end{align*}
Thus, we have $f\equiv m$ on $B(x_0,2r)$.
Repeating this process, we obtain $f\equiv m$ on $B(x_0,kr)$
for any $k\in\mathbb{N}$, which further implies that
$f\equiv m$ on $\mathbb{R}^n$.
On the other hand, since $f\in C_{\mathrm{c}}^2$,
it follows that
there exists $z\in\mathbb{R}^n$ such that $f(z)=0$,
and hence $f\equiv m=0$,
which contradicts $f\not\equiv 0$. Therefore, \eqref{equ-f-fav}
and \eqref{equ-f-B} hold. This completes the proof of
Theorem \ref{wenti1fouding}.
\end{proof}

\subsection{Proof of Theorem \ref{ruofanshu1}}\label{sub-weak-min}
In this subsection, we prove Theorem \ref{ruofanshu1}. To this end,
we introduce the following notation.
For any Lebesgue measurable subset $H$ of
$\mathbb{R}^{2n}=\mathbb{R}^n\times \mathbb{R}^n$ and
for any $r\in(0,\infty)$,
let
\begin{align}\label{nur}
\nu_{r}(H):=
\iint_{(x,y)\in H,x\ne y}
\rho_r(x-y)\,dx\,dy.
\end{align}
For any $f\in L_{\rm{loc}}^1$, $r\in(0,\infty)$, and $\lambda\in(0,\infty)$,
let
\begin{align}\label{Hru}
H_{\lambda,r}(f):=\left\{
(x,y)\in\mathbb{R}^n\times\mathbb{R}^n
: x\ne y, \frac{|f(y)-f_{B(x,r)}|}{|x-y|}>\lambda
\right\}.
\end{align}
Now, we prove Theorem \ref{ruofanshu1}.
\begin{proof}[Proof of Theorem \ref{ruofanshu1}]
Choose $f\in C_{\mathrm{c}}^\infty(\mathbb{R})$
such that
$f(x):=x^2$ for any $x\in [0,1]$.
For any $r\in(0,\infty)$,
let $\rho_r$ be the same as in \eqref{teshu} with $\theta\in(0,1)$.
Fix $r\in(0,\frac{1}{4})$, $x\in (\frac{1}{4},\frac{1}{2})$,
and $y\in (x,\sqrt{x^2+\frac{r^2}{6}})$. Then we have
\begin{align*}
|y-x|\le \sqrt{x^2+\frac{r^2}{6}}-x
\le \frac{r}{\sqrt{6}}<r.
\end{align*}
By this, we obtain
\begin{align}\label{equ-nu-H}
\nu _r(H_{\lambda,r}(f))\ge
\int_{\frac{1}{4}}^{\frac{1}{2}}
{\int_{\{y\in(x,\sqrt{x^2+\frac{r^2}{6}}):
\frac{|f(y)-f_{B(x,r)}|}
{|x-y|}>\lambda\}}{\frac{\theta}
{2r^{\theta}|x-y|^{1-\theta}}\,dy\,dx}}.
\end{align}
Since $x\in (\frac{1}{4},\frac{1}{2})$
and $r\in(0,\frac{1}{4})$,
it follows that
\begin{align}\label{equ-r2}
\sqrt{x^2+\frac{r^2}{6}}-x
=\frac{r^2}{6(\sqrt{x^2+\frac{r^2}{6}}+x)}
\geq\frac{r^2}{18}.
\end{align}
On the other hand, a direct calculation shows that
$$|f(y)-f_{B(x,r)}|
=\left|y^2-x^2-\frac{r^2}{3}\right|
\geq\frac{r^2}{6}.$$
From this, \eqref{equ-nu-H}, \eqref{equ-r2},
and a change of variables, we deduce that,
for any $\lambda\in(3,\infty)$,
\begin{align*}
\nu _r(H_{\lambda,r}(f))&\ge \int_{\frac{1}{4}}^{\frac{1}{2}}
{\int_{\{ y\in(x,\sqrt{x^2+\frac{r^2}{6}}):
\frac{r^2}{6|x-y|}>\lambda\}}
{\frac{\theta}{2r^{\theta}|x-y|^{1-\theta}}\,dy\,dx}}\\
&\gtrsim
\int_{\{ y\in(0,\min\{\frac{r^2}{18},
\frac{r^2}{6\lambda}\}) \}}
\frac{1 }{r^{\theta }y^{1-\theta}}\,dy
\sim\lambda^{-\theta}r^{\theta}.
\end{align*}
Since $\theta<1\le q$, it follows that
\begin{align*}
\sup_{\lambda\in(0,\infty)}
\lambda\left[\nu_r(H_{\lambda,r}(f))\right]^\frac{1}{q}
\gtrsim
\sup_{\lambda\in(3,\infty)}
\lambda^{1-\frac{\theta}{q}}r^{\frac{\theta}{q}}
=\infty.
\end{align*}
This completes the proof of Theorem \ref{ruofanshu1}.
\end{proof}
\section{Proof of Theorem \ref{tailek}}\label{s5}
In this section, we prove Theorem \ref{tailek}.
To this end, we recall the concept of
Taylor polynomials.
Let $k\in\mathbb{N}$ and $f\in C^k$.
Then, by the Taylor remainder
theorem, we have, for any $x,h\in\mathbb{R}^n$,
\begin{align}\label{tailezhankaibiaozhun}
f(x+h)=\sum_{\{\beta\in\mathbb{Z}^n_+:|\beta|\leq k\}}
\frac{\partial^\beta f(x)}{\beta!}h^\beta+R^{(k)}(x,h),
\end{align}
where $R^{(k)}(x,h)$ is the
\textit{k-th Taylor remainder} of $f$.

By the definitions of $\sigma_k$ and $D^k(f)$,
we rewrite \eqref{tailezhankaibiaozhun}
in a vectorized form as follows
(see, for instance, \cite[Section 2]{fkr}):
\begin{align}\label{tailelingyizhong}
f(x+h)=f(x)+\sum_{m=1}^{k}\frac{1}{m!}
\sigma _m(h)\nabla\left(D^{m-1}f\right)(x)h+R^{(k)}(x,h),
\end{align}
where the remainder $R^{(k)}(x,h)$ is given by
\begin{align*}
R^{(k)}(x,h):=\sum_{\{\beta\in\mathbb{Z}_+^n:|\beta|=k\}}
\frac{\partial^{\beta}f(\widetilde{y})
-\partial^{\beta}f(x)}{\beta!}h^{\beta}
\end{align*}
with some $\widetilde{y}$ on the line segment joining
$x$ and $x+h$
(see, for instance, \cite[Section 2]{fkr}).
Then it is easy to show that, for
any $x\in\mathbb{R}^n$,
\begin{equation}\label{taileduo}
T_x^{k}(f)(\cdot):=f(x)+
\sum_{m=1}^{k}\frac{1}{m!}\sigma_m(\cdot-x)
\nabla\left(D^{m-1}f\right)(x)(\cdot-x)\in\mathcal{P}_k.
\end{equation}

The following concept of the
associate space of a ball Banach function space can be found in
\cite[p.9]{shyy}.
\begin{definition}
Let $X$ be a ball Banach function space. The \textit{associate space
$X'$} of $X$ is defined by setting
$X':=\{f\in\mathscr{M}:\|f\|_{X'}<\infty\}$,
where, for any $f\in X'$,
$$\|f\|_{X'}:=\sup\left\{\|fg\|_{L^1}:\,\|g\|_X=1\right\}$$
and $\|\cdot\|_{X'}$ is called the associate norm of $\|\cdot\|_X$.
\end{definition}
We also recall the definition of the absolutely continuous norm.
\begin{definition}
A ball Banach function space $X$ is said to have an \textit{absolutely
continuous norm} if, for any $f\in X$ and any sequence of measurable
sets $\{E_j\}_j$ such that $\mathbf{1}_{E_j}\to0$ a.e. as
$j\to\infty$, we have $\|f\mathbf{1}_{E_j}\|_X\to0$ as $j\to\infty$.
\end{definition}
To prove Theorem \ref{tailek}, we need the
following two lemmas.
The first lemma is precisely \cite[Proposition 3.18]{jyy25}.
\begin{lemma}\label{yiyou}
Suppose that the assumptions of Theorem \ref{tailek} hold.
Then $f\in W^{k,X}$ if and only
if $f\in W^{k-1,X}$ and
there exists a function $H$ satisfying $H(x)\in\mathbb{M}_{n^{k-1},n}$
for almost every $x\in\mathbb{R}^n$ and
$|H(\cdot)|\in X$ such that
\begin{align*}
\lim_{r\rightarrow0^+}\left\|\frac{1}{r^k}
\left[\fint_{B(\cdot,r)}\left|f(y)-T_{\cdot}^{k-1}(f)(y)
-\frac{1}{k!}\sigma_k(y-\cdot)H(\cdot)(y-\cdot)\right|^q
\,dy\right]^{\frac{1}{q}}\right\|_X=0.
\end{align*}
\end{lemma}
We also need the following lemma, which is a part of \cite[Theorem 3]{bh93}.
\begin{lemma}\label{zybds1}
Let $m\in\mathbb{N}$ and $f\in W^{m,L^1}_{\rm{loc}}$.
Then there exists a positive constant $C$, depending only on $m$
and $n$, such that,
for any $\beta\in\{\alpha\in\mathbb{Z}_+^n:|\alpha|\le m-1\}$
and for almost every $x,y\in\mathbb{R}^n$ satisfying
$|f(y)|<\infty$ and $|\partial^\beta f(x)|<\infty$,
\begin{equation*}
\left|f(y)-T_x^{m-1}(f)(y)\right|\leq
C\left[\mathcal{M}_{|x-y|}(|\nabla^mf|)(y)
+\mathcal{M}_{|x-y|}(|\nabla^mf|)(x)\right]|x-y|^m,
\end{equation*}
where, for
any $R\in(0,\infty)$, $g\in L_{\mathrm{loc}}^1$,
and $x\in\mathbb{R}^n$,
$\mathcal{M}_{R}(g)(x):=\sup_{r\in(0,R)}\fint_{B(x,r)}g(y)
\,dy$.
\end{lemma}
We also establish the following
higher-order BBM-type formula involving Taylor polynomials,
which plays an essential role in the proof of Theorem \ref{tailek}.
\begin{theorem}\label{taile}
Let $q\in[1,\infty)$, $k\in\mathbb{N}$, and
$X$ be a ball Banach function space.
Let $\{\rho_r\}_{r\in(0,\infty)}$ satisfy
the conditions of Definition \ref{djh}.
Assume that, for any given $M\in(0,\infty)$,
there exists $N\in(2M+1,\infty)$ such that
\begin{align}\label{rou}
\lim_{r\rightarrow0^+}
\left\|\left[\int_{B(\mathbf{0},M)}
\frac{\rho_r(\cdot-y)}{|\cdot-y|^{kq}}
\,dy\right]^{\frac{1}{q}}
\mathbf{1}_{\mathbb{R}^n\setminus B( 0,N )}
\right\|_{X}=0.
\end{align}
Then, for any $f\in C_{\mathrm{c}}^{\infty}$,
\begin{align}\label{equ-Irr}
&\lim_{r\to 0^+}
\left\|
\left[\int_{\mathbb{R}^n}
\frac{|f(y)-T_{\cdot}^{k-1}(f)(y)|^q}{|\cdot-y|^{kq}}
\rho_r(\cdot-y)\,dy\right]^{\frac{1}{q}}
\right\|_{X}\notag\\
&\quad=\frac{1}{k![\sigma(\mathbb{S}^{n-1})]^{\frac{1}{q}}}
\left\|\left[\int_{\mathbb{S}^{n-1}}\left|\sigma_k(\xi)
\nabla\left(D^{k-1}f\right)(\cdot)\xi\right|^qd\sigma(\xi)
\right]^{\frac{1}{q}}
\right\|_{X}.
\end{align}
\end{theorem}
\begin{proof}
Fix $k\in\mathbb{N}$ and $f\in C_{\mathrm{c}}^{\infty}$,
and let $\{\rho_r\}_{r\in(0,\infty)}$
be a family of functions
satisfying \eqref{rou} and the conditions of Definition \ref{djh}.
Since $f\in C^k_{\rm{c}}$,
it follows that
there exists $M\in(0,\infty)$ such that, for any
$m\in\{0,\ldots,k\}$,
$\mathop{\mathrm{supp}}(|\nabla^mf|)\subset B(\mathbf{0},M)$.
Fix $r\in(0,\infty)$, and let
\begin{align*}
\text{I}_r:=\left\|\left[\int_{\mathbb{R}^n}
\frac{|f(y)-T_{\cdot}^{k-1}(f)(y)|^q}{|\cdot-y|^{kq}}
\rho_r(\cdot-y)\,dy\right]^{\frac{1}{q}}\right\|_{X}.
\end{align*}
Fix $\varepsilon\in(0,\infty)$. Then
there exists $\delta\in(0,\min\{\varepsilon,1\})$ such that,
for any $x,y\in\mathbb{R}^n$
with $|x-y|<\delta$,
\begin{align}\label{lianxukewei}
\sum_{\{\beta\in\mathbb{Z}_+^n:|\beta|=k\}}
\left| \frac{\partial^{\beta}f(y)-\partial
^{\beta}f(x)}{\beta!}\right|<\varepsilon.
\end{align}
Choose $N\in (2M+1,\infty)$.
By the assumption that $q\in[1,\infty)$, we can decompose
$\mathrm{I}_r$ as follows:
\begin{align*}
\text{I}_r&\leq\left\|\left[\int_{|\cdot-y|<\delta}
\frac{|f(y)-T_{\cdot}^{k-1}(f)(y)|^q}{|\cdot-y|^{kq}}
\rho_r(\cdot-y)\,dy\right]^{\frac{1}{q}}
\right\|_{X}\\
&\quad+\left\|\left[\int_{\delta\leq|\cdot-y|<2N}
\cdots\right]^{\frac{1}{q}}\right\|_{X}+
\left\|\left[\int_{|\cdot-y|\geq2N}\cdots\right]^{\frac{1}{q}}
\right\|_{X}\\
&=:\text{J}_1+\text{J}_2+\text{J}_3.
\end{align*}

We first estimate $\text{J}_1$. Fix $x\in\mathbb{R}^n$.
By \eqref{tailelingyizhong} and Lemma
\ref{fmsbds}, we conclude that, for any $z\in(0,1)$,
\begin{align*}
&\int_{|x-y|<\delta}
\frac{|f(y)-T_{x}^{k-1}(f)(y)|^q}{|x-y|^{kq}}
\rho_r(x-y)\,dy\\
&\quad=\int_{|x-y|<\delta}
\frac{|\frac{1}{k!}\sigma_k(y-x)\nabla
(D^{k-1}f)(x)(y-x)
+R_{f}^{(k)}(x,y-x)|^q}{|x-y|^{kq}}
\rho_r(x-y)\,dy\\
&\quad\leq(1+z)^q\int_{|x-y|<\delta}
\frac{|\frac{1}{k!}\sigma_k(y-x)\nabla
(D^{k-1}f)(x)(y-x)|^q}{|x-y|^{kq}}
\rho_r(x-y)\,dy\\
&\quad\quad+\left(1+\frac{1}{z}\right)^q\int_{|x-y|<\delta }
\frac{|R_{f}^{(k)}(x,y-x)|^q}{|x-y|^{kq}}
\rho_r(x-y)\,dy\\
&\quad=:\left(\frac{1+z}{k!}\right)^qM_1+\left(1+\frac{1}{z}\right)^qM_2.
\end{align*}
First, we estimate $M_1$. A change of variables implies that
\begin{align*}
M_1&=\int_{B(\mathbf{0},\delta)}
\frac{|\sigma _k(h)\nabla(D^{k-1}f)(x)h|^q}{|h|^{kq}}
\rho_r(h)\,dh\\
&=\int_{\mathbb{S}^{n-1}}
\int_0^{\delta}
\frac{|\sigma_k(\xi)\nabla(D^{k-1}f)(x)\xi|^q{\eta}^{kq}}
{|\eta\xi|^{kq}}
\rho_r(\eta\xi){\eta}^{n-1}\,d\eta\,d\sigma(\xi)\\
&=\int_{\mathbb{S}^{n-1}}\int_0^{\delta}
\left|\sigma_k(\xi)\nabla\left(D^{k-1}f\right)(x)
\xi\right|^q
\rho_r(\eta\xi){\eta}^{n-1}\,d\eta\,d\sigma(\xi).
\end{align*}
Choose $\xi_1\in\mathbb{S}^{n-1}$.
Since $\rho_r$ is radial,
it follows that $\rho_r(\eta\xi)=\rho_r(\eta\xi_1)$
for any $\xi\in\mathbb{S}^{n-1}$.
This, together with the assumption
that $\int_{\mathbb{R}^n}\rho_r(x)\,dx=1$,
further implies that
\begin{align}\label{jingxiang}
M_1&=\int_{\mathbb{S}^{n-1}}\int_0^{\delta}
\left|\sigma_k(\xi)\nabla\left(D^{k-1}f\right)(x)
\xi\right|^q
\rho_r(\eta\xi_1){\eta}^{n-1}
\,d\eta\,d\sigma(\xi)\notag\\
&\leq\int_{\mathbb{S}^{n-1}}
\left|\sigma_k(\xi)\nabla\left(D^{k-1}f\right)(x)
\xi\right|^q\,d\sigma(\xi)
\int_0^{\infty}\rho_r(\eta\xi_1){\eta}^{n-1}\,d\eta\notag\\
&=\frac{1}{\sigma(\mathbb{S}^{n-1})}
\int_{\mathbb{S}^{n-1}}
\left|\sigma_k(\xi)\nabla
\left(D^{k-1}f\right)(x)\xi\right|^q
\,d\sigma(\xi)
\int_{\mathbb{S}^{n-1}}\int_0^{\infty}
\rho_r(\eta\zeta){\eta}^{n-1}\,d\eta\,
d\sigma(\zeta)\notag\\
&\leq\frac{1}{\sigma(\mathbb{S}^{n-1})}
\int_{\mathbb{S}^{n-1}}\left|\sigma_k(\xi)
\nabla\left(D^{k-1}f\right)(x)\xi\right|^q\,
d\sigma(\xi).
\end{align}
Next, we estimate
$M_2$. By the definition of $R_f^{(k)}$ and by
\eqref{lianxukewei},
we conclude that, for any $y\in B(x,\delta)$,
$$\left|R_f^{(k)}(x,y-x)\right|
=\left|\sum_{\{\beta\in\mathbb{Z}_+^n:|\beta|=k\}}
\frac{\partial^{\beta}f(\widetilde{y})-
\partial^{\beta}f(x)}{\beta!}
(y-x)^\beta\right|\leq\varepsilon|y-x|^k,$$
where $\widetilde{y}$ is some point in $B(x,\delta)$.
Since $\int_{\mathbb{R}^n}\rho_{r}(x)\,dx=1$, it follows
that
\begin{align}\label{equ-Remain}
\int_{|x-y|<\delta }
\frac{|R_{f}^{(k)}(x,y-x)|^q}{|x-y|^{kq}}
\rho_r(x-y)\,dy\leq
\int_{|x-y|<\delta}{\varepsilon}^q\rho_r(x-y)\,
dy\leq{\varepsilon}^q.
\end{align}
This, combined with \eqref{basic-inequality} and
\eqref{jingxiang},
further implies that, for any $z\in(0,1)$,
\begin{align*}
&\left[\int_{|x-y|<\delta}
\frac{|f(y)-T_{x}^{k-1}(f)(y)|^q}{|x-y|^{kq}}
\rho_r(x-y)\,dy\right]^{\frac{1}{q}}\\
&\quad\leq\left[\frac{(1+z)^q}{(k!)^q\sigma
(\mathbb{S}^{n-1})}\int_{\mathbb{S}^{n-1}}
\left|\sigma _k(\xi)\nabla(D^{k-1}f)(x)\xi\right|
^qd\sigma(\xi)+\left(1+\frac{1}{z}\right)^q{\varepsilon}^q
\right]^{\frac{1}{q}}\\
&\quad\leq\frac{1+z}{k![\sigma (\mathbb{S}^{n-1})]
^{\frac{1}{q}}}\left[\int_{\mathbb{S}^{n-1}}
\left|\sigma_k(\xi)\nabla(D^{k-1}f)(x)\xi\right|
^qd\sigma(\xi)\right]^{\frac{1}{q}}+
\left(1+\frac{1}{z}\right)\varepsilon.
\end{align*}
Note that, if $|x|>M+1$,
then, for any $y\in B(x,\delta)$,
we have $|y|>M$, which, together with the fact that
$\mathop{\mathrm{supp}}(|\nabla^mf|)\subset B(\mathbf{0},M)$
for any $m\in\{0,\ldots,k\}$,
further implies that $|f(y)-T_x^{k-1}(f)(y)|=0$.
By this and the definition of ball Banach function spaces,
we conclude that,
for any $z\in(0,1)$,
\begin{align}\label{j1s}
\text{J}_1\leq\frac{1+z}{k![\sigma
(\mathbb{S}^{n-1})]^{\frac{1}{q}}}\left\|
\left[\int_{\mathbb{S}^{n-1}}\left|\sigma_k(\xi)
\nabla(D^{k-1}f)(\cdot)\xi\right|^qd\sigma(\xi)
\right]^{\frac{1}{q}}
\right\|_{X}
+\left(1+\frac{1}{z}\right)\varepsilon
\left\|\mathbf{1}_{B(\mathbf{0},M+1)}\right\|_{X}.
\end{align}
Next, we establish a lower bound for $\rm{J}_1$. To this end,
fix $x\in\mathbb{R}^n$. Using Lemma \ref{fmsbds}, we conclude that, for any
$a,b\in\mathbb{C}$ and $z\in(0,1)$,
$$(1+z)^q||a|-|b||^q+\frac{(1+z)^q}{z^q}|b|^q
\geq||a|-|b|+|b||^q,$$
which further implies that
\begin{align*}
|a+b|^q\geq||a|-|b||^q\geq\frac{1}{(1+z)^q}
|a|^q-\frac{1}{z^q}|b|^q.
\end{align*}
By this, we obtain
\begin{align}\label{equ-Bxdelta}
&\int_{|x-y|<\delta}
\frac{|\frac{1}{k!}\sigma_k(y-x)
\nabla(D^{k-1}f)(x)(y-x)
+R_{f}^{(k)}(x,y-x)|^q}{|x-y|^{kq}}
\rho_r(x-y)\,dy\notag\\
&\quad\geq\frac{1}{(1+z)^q}\int_{|x-y|<\delta}
\frac{|\frac{1}{k!}\sigma_k(y-x)
\nabla(D^{k-1}f)(x)(y-x)|^q}{|x-y|^{kq}}
\rho_r(x-y)\,dy\notag\\
&\quad\quad-\frac{1}{z^q}\int_{|x-y|<\delta }
\frac{|R_{f}^{(k)}(x,y-x)|^q}{|x-y|^{kq}}
\rho_r(x-y)\,dy.
\end{align}
First, we have
\begin{align}\label{equ-H1H2}
&\frac{1}{(1+z)^q}\int_{|x-y|<\delta }
\frac{|\frac{1}{k!}\sigma_k(y-x)
\nabla(D^{k-1}f)(x)(y-x)|^q}{|x-y|^{kq}}
\rho_r(x-y)\,dy\notag\\
&\quad=\frac{1}{(1+z)^q(k!)^q}\left[\int_{\mathbb{R}^n}
\frac{|\sigma_k(h)\nabla(D^{k-1}f)(x)h|^q}{|h|^{kq}}
\rho_r(h)\,dh-\int_{\mathbb{R}^n\setminus B(\mathbf{0},\delta )}
\cdots\right]\notag\\
&\quad=:\frac{1}{(1+z)^q(k!)^q}(H_1-H_2).
\end{align}
Applying
the assumption that $\int_{\mathbb{R}^n}\rho_r(y)\,dy=1$
and an argument similar to that used in \eqref{jingxiang}, we obtain
\begin{align}\label{equ-H1}
H_1&=\int_{\mathbb{S}^{n-1}}\int_0^{\infty}\left|\sigma_k(\xi)
\nabla\left(D^{k-1}f\right)(x)\xi\right|^q
\rho_r(\eta\xi_1){\eta}^{n-1}\,d\eta\,d\sigma(\xi)\notag\\
&=\int_{\mathbb{S}^{n-1}}\left|\sigma_k(\xi)
\nabla\left(D^{k-1}f\right)(x)\xi\right|^q\,d\sigma(\xi)
\int_0^{\infty}\rho_r(\eta\xi_1){\eta}^{n-1}\,d\eta\notag\\
&=\frac{1}{\sigma(\mathbb{S}^{n-1})}
\int_{\mathbb{S}^{n-1}}\left|\sigma_k(\xi)
\nabla\left(D^{k-1}f\right)(x)\xi\right|^qd\sigma(\xi)
\int_{\mathbb{S}^{n-1}}\int_0^{\infty}\rho_r(\eta\zeta)
{\eta}^{n-1}\,d\eta\,d\sigma(\zeta)\notag\\
&=\frac{1}{\sigma(\mathbb{S}^{n-1})}
\int_{\mathbb{S}^{n-1}}\left|\sigma_k(\xi)
\nabla\left(D^{k-1}f\right)(x)\xi\right|^q\,d\sigma(\xi),
\end{align}
where $\xi_1$ is some point in $\mathbb{S}^{n-1}$.
For $H_2$, by \eqref{tailelingyizhong}, we obtain
\begin{align*}
H_2&\leq(k!)^q\int_{\mathbb{R}^n
\setminus B(\mathbf{0},\delta )}
\frac{[\sum_{\{\beta\in\mathbb{Z}^n_+:|\beta|=k\}}
\frac{|{\partial}^{\beta}f(x)}
{\beta!}|h|^k]^q}{|h|^{kq}}
\rho_r(h)\,dh\\
&=(k!)^q\left[\sum_{\{\beta\in\mathbb{Z}^n_+:|\beta|=k\}}
\frac{|{\partial}^{\beta }f(x)|}
{\beta!}\right]^q\int_{\mathbb{R}^n\setminus
B(\mathbf{0},\delta)}\rho_r(h)\,dh.
\end{align*}
From this, \eqref{equ-Remain}, \eqref{equ-Bxdelta},
\eqref{equ-H1H2},
and \eqref{equ-H1},
we infer that
\begin{align*}
&\int_{|x-y|<\delta }
\frac{|\frac{1}{k!}\sigma _k(y-x)
\nabla(D^{k-1}f)(x)(y-x)
+R_{f}^{(k)}(x,y-x)|^q}{|x-y|^{kq}}
\rho_r(x-y)\,dy\\
&\quad\geq\frac{1}{(1+z)^q(k!)^q\sigma(\mathbb{S}^{n-1})}
\int_{\mathbb{S}^{n-1}}\left|\sigma _k(\xi )
\nabla\left(D^{k-1}f\right)(x)\xi\right|^q\,d\sigma (\xi )\\
&\quad\quad -\left[\sum_{\{\beta\in\mathbb{Z}^n_+:|\beta|=k\}}\frac{|{\partial}^{\beta }f(x)|}
{\beta!(1+z)}\right]^q\int_{\mathbb{R}^n\setminus
B(\mathbf{0},\delta )}
\rho_r(h)\,dh-\frac{\varepsilon ^q}{z^q}.
\end{align*}
This, together with the fact that if $|x|>M+1$,
then $\nabla(D^{k-1}f)(x)=0$, further implies that
\begin{align*}
\text{J}_1
&\geq\frac{1}{(1+z)k![\sigma (\mathbb{S}^{n-1})]
^{\frac{1}{q}}}\left\|\left[\int_{\mathbb{S}^{n-1}}
\left|\sigma _k(\xi)\nabla\left(D^{k-1}f\right)(\cdot)\xi
\right|^qd\sigma(\xi)\right]^{\frac{1}{q}}
\right\|_{X}\\
&\quad-\frac{\varepsilon }{z}
\left\|\mathbf{1}_{B(\mathbf{0},M+1)}\right\|
_{X}-\frac{1}{1+z}\left[\int_{\mathbb{R}^n\setminus
B(\mathbf{0},\delta)}
\rho_r(h)\,dh\right]^{\frac{1}{q}}\left\|\sum_{|\beta|=k }
\frac{|{\partial}^{\beta}f(\cdot)|}{\beta !}
\mathbf{1}_{B(\mathbf{0},M+1)}\right\|_{X}.
\end{align*}
From this, the assumption that
$\lim_{r\to 0^+}
\int_{\mathbb{R}^n\setminus B(\mathbf{0},\delta)}\rho_r(h)\,dh=0$, and
the fact that $|\nabla^kf|\in C_{\mathrm{c}}\subset X$,
we deduce that
\begin{align}\label{equ-lower-J1}
&\liminf_{r\to 0^+}\mathrm{J}_1+\frac{\varepsilon}{z}
\left\|\mathbf{1}_{B(\mathbf{0},M+1)}\right\|_{X}\notag\\
&\quad\ge\frac{1}{(1+z)k![\sigma (\mathbb{S}^{n-1})]
^{\frac{1}{q}}}\left\|\left[\int_{\mathbb{S}^{n-1}}
\left|\sigma_k(\xi)\nabla\left(D^{k-1}f\right)(\cdot) \xi
\right|^q\,d\sigma(\xi)\right]^{\frac{1}{q}}
\right\|_{X}.
\end{align}

We next estimate $\text{J}_2$.
For any $x\in \mathbb{R}^n \setminus B(\mathbf{0},M+2N)$
and $y\in B(x,2N)\setminus B(x,\delta)$,
we have $y\in \mathbb{R}^n \setminus B(\mathbf{0},M)$,
which, combined with the fact that $\mathop{\mathrm{supp}}(f)\subset B(\mathbf{0},M)$,
further implies that
$|f(y)-T_{x}^{k-1}(f)(y)|=0$.
By this, we obtain
\begin{align}\label{equ-J2}
\text{J}_2&=\left\|\left[\int_{\delta \leq |\cdot-y|<2N }
\frac{|f(y)-T_{\cdot}^{k-1}(f)(y)|^q}{|\cdot-y|^{kq}}
\rho_r(\cdot-y)\,dy\right]^{\frac{1}{q}}
\mathbf{1}_{B(\mathbf{0},M+2N)}
\right\|_{X}.
\end{align}
From the Taylor remainder theorem and \eqref{tailelingyizhong},
we infer that,
for any $x$ and $y\in\mathbb{R}^n$,
\begin{align}\label{wucha}
\left|f(y)-T_{x}^{k-1}(f)(y)\right|
&=\left|\sum_{\{\beta\in\mathbb{Z}^n_+:|\beta|=k\}}
\frac{\partial ^{\beta}f(x)}
{\beta!}(y-x)^\beta+R_f^{(k)}(x,y-x)
\right|\notag\\
&\lesssim\sum_{\{\beta\in\mathbb{Z}^n_+:|\beta|=k\}}
\frac{\|\partial^{\beta}f\|_{L^\infty}}{\beta!}
|x-y|^k\lesssim|x-y|^k
\end{align}
with the implicit positive constants independent of $x$ and $y$.
Using this, \eqref{equ-J2}, the fact
that $\mathbf{1}_{B(\mathbf{0},M+2N)}\in X$
and $\lim_{r\to 0^+}\int_{B(\mathbf{0},\delta)^\complement}
\rho_r(y)\,dy=0$,
and a change of variables, we obtain
\begin{align}\label{j2}
\limsup_{r\to 0^+}\text{J}_2&\lesssim \limsup_{r\to 0^+}
\left\|\left[\int_{\delta
\leq |\cdot-y|<2N }\rho_r(\cdot-y)\,dy\right]
^{\frac{1}{q}}\mathbf{1}_{B(\mathbf{0},M+2N)}\right\|_{X}\notag\\
&\leq\limsup_{r\to 0^+}\left[\int_{B(\mathbf{0},\delta)^\complement}
\rho_r(y)\,dy\right]
^{\frac{1}{q}}\|\mathbf{1}_{B(\mathbf{0},M+2N)}\|_X=0.
\end{align}

Finally, we estimate $\text{J}_3$. By the assumption that
$\mathop{\mathrm{supp}}(f)\subset B(\mathbf{0},M)$, we find that, for
any $m\in\{1,\ldots,k-1\}$ and $x\in\mathbb{R}^n\setminus B(\mathbf{0},M)$,
$$\left|\nabla\left(D^{m-1}f\right)(x)\right|=0.$$
Thus, we can decompose $\rm{J}_3$ as follows:
\begin{align*}
\text{J}_3&\leq\left\|\left[\int_{|\cdot-y|\geq 2N }
\frac{|f(y)-T_{\cdot}^{k-1}(f)(y)|^q}{|\cdot-y|^{kq}}
\rho_r(\cdot-y)\,dy\right]^{\frac{1}{q}}
\mathbf{1}_{B(\mathbf{0},M)}
\right\|_{X}\\
&\quad +\left\|\left[\int_{|\cdot-y|\geq 2N }
\frac{|f(y)|^q}{|\cdot-y|^{kq}}
\rho_r(\cdot-y)\,dy\right]^{\frac{1}{q}}
\mathbf{1}_{B(\mathbf{0},N)\setminus B(\mathbf{0},M)}
\right\|_{X}\\
&\quad +\left\|\left[\int_{|\cdot-y|\geq 2N }
\frac{|f(y)|^q}{|\cdot-y|^{kq}}
\rho_r(\cdot-y)\,dy\right]^{\frac{1}{q}}
\mathbf{1}_{\mathbb{R}^n\setminus B(\mathbf{0},N)}
\right\|_{X}\\
&=: \text{T}_1+\text{T}_2+\text{T}_3.
\end{align*}
We first estimate $T_1$. By \eqref{wucha}, the assumption on $\rho_r$,
and a change of variables, we obtain
\begin{align}\label{equ-T1}
\limsup_{r\to 0^+}\text{T}_1&=\limsup_{r\to 0^+}
\left\|\left[\int_{|\cdot-y|\geq 2N }
\frac{|f(y)-T_{\cdot}^{k-1}(f)(y)|^q}{|\cdot-y|^{kq}}
\rho_r(\cdot-y)\,dy\right]^{\frac{1}{q}}
\mathbf{1}_{B(\mathbf{0},M)}
\right\|_{X}\notag\\
&\lesssim\limsup_{r\to 0^+}\left[\int_{|h|\geq 2N }
\rho_r(h)\,dh\right]^{\frac{1}{q}}
\left\|
\mathbf{1}_{B(\mathbf{0},M)}
\right\|_{X}=0.
\end{align}
For $\text{T}_2$, using a change of variables, we have
\begin{align}\label{T3}
\limsup_{r\to 0^+}\text{T}_2&
\leq \limsup_{r\to 0^+}\frac{\|f\|_{L^\infty}}{(2N)^k}
\left\|\left[\int_{|\cdot-y|\geq2N }
\rho_r(\cdot-y)\,dy\right]^{\frac{1}{q}}
\mathbf{1}_{B(\mathbf{0},N)\setminus B(\mathbf{0},M)}
\right\|_{X}\notag\\
&=\limsup_{r\to 0^+}\frac{\|f\|_{L^\infty(\mathbb{R}^n)}}
{(2N)^{k}}\left[\int_{|y|\geq2N }
\rho_r(y)\,dy\right]^{\frac{1}{q}}
\left\|\mathbf{1}_{B(\mathbf{0},N)\setminus B(\mathbf{0},M)}
\right\|_{X}=0.
\end{align}
For $\text{T}_3$, from the fact that $\mathop{\mathrm{supp}}(f)
\subset B(\mathbf{0},M)$ and \eqref{rou}, we deduce that
\begin{align*}
\limsup_{r\rightarrow 0^+}\text{T}_3&
\leq\limsup_{r\rightarrow 0^+}
\lVert f \rVert _{L^{\infty}( \mathbb{R}^n )}
\left\|\left[\int_{B(\mathbf{0},M)}
\frac{\rho_r(\cdot-y)}{|\cdot-y|^{kq}}
\,dy\right]^{\frac{1}{q}}
\mathbf{1}_{\mathbb{R}^n\setminus B(\mathbf{0},N)}
\right\|_{X}=0.
\end{align*}
This, together with \eqref{equ-T1} and \eqref{T3}, further implies that
$\limsup_{r\to 0^+}\mathrm{J}_3=0$.
By this, \eqref{j1s}, and \eqref{j2},
we find that
\begin{align}\label{equ-limsup-Ir}
\limsup_{r\to0^+}
\mathrm{I}_r&\le
\frac{1+z}{k![\sigma
(\mathbb{S}^{n-1})]^{\frac{1}{q}}}\left\|
\left[\int_{\mathbb{S}^{n-1}}\left|\sigma _k(\xi )
\nabla ( D^{k-1}f ) (\cdot  ) \xi\right|^qd\sigma (\xi )
\right]^{\frac{1}{q}}
\right\|_{X}\notag\\
&\quad+\left(1+\frac{1}{z}\right)\varepsilon
\left\|\mathbf{1}_{B(\mathbf{0},M+1)}\right\|_{X}.
\end{align}
Letting $\varepsilon\to0^+$ first
and then $z\to0^+$ in \eqref{equ-limsup-Ir}, we obtain
\begin{align}\label{equ-upper-Ir}
\limsup_{r\to0^+}
\mathrm{I}_r
\le\frac{1}{k![\sigma
(\mathbb{S}^{n-1})]^{\frac{1}{q}}}\left\|
\left[\int_{\mathbb{S}^{n-1}}\left|\sigma _k(\xi )
\nabla(D^{k-1}f)(\cdot)\xi\right|^qd\sigma (\xi )
\right]^{\frac{1}{q}}
\right\|_{X}.
\end{align}
On the other hand, from \eqref{equ-lower-J1}, we infer that
\begin{align*}
\liminf_{r\to 0^+}\mathrm{I}_r
\ge\frac{1}{(1+z)k![\sigma (\mathbb{S}^{n-1})]
^{\frac{1}{q}}}\left\|\left[\int_{\mathbb{S}^{n-1}}
\left|\sigma _k(\xi)\nabla\left(D^{k-1}f\right)(\cdot)\xi
\right|^q\,d\sigma(\xi)\right]^{\frac{1}{q}}
\right\|_{X}-\frac{\varepsilon }{z}
\left\|\mathbf{1}_{B(\mathbf{0},M+1)}\right\|_{X}.
\end{align*}
Letting $\varepsilon\to0^+$ first and then $z\to0^+$, we conclude that
\begin{align*}
\liminf_{r\to 0^+}\mathrm{I}_r
\ge\frac{1}{k![\sigma(\mathbb{S}^{n-1})]
^{\frac{1}{q}}}\left\|\left[\int_{\mathbb{S}^{n-1}}
\left|\sigma _k(\xi)\nabla\left(D^{k-1}f\right)(\cdot)\xi
\right|^q\,d\sigma(\xi)\right]^{\frac{1}{q}}
\right\|_{X}.
\end{align*}
This, combined with \eqref{equ-upper-Ir},
further implies \eqref{equ-Irr},
which completes the proof of Theorem \ref{taile}.
\end{proof}
We now prove Theorem \ref{tailek}.
\begin{proof}[Proof of Theorem \ref{tailek}]
We first show the sufficiency (``if'' part).
Fix $f\in W^{k-1,X}$
such that, for almost
every $x\in\mathbb{R}^n$, there exists $H(x)\in
\mathbb{M}_{n^{k-1},n}$ with $|H|\in X$ such
that \eqref{equ-f-taylor} holds. By this, we obtain
\begin{align*}
&\limsup_{r\rightarrow0^+}\left\|\frac{1}{r^k}
\left[\fint_{B(\cdot,r)}\left|f(y)-T_{\cdot}^{k-1}(f)(y)
-\frac{1}{k!}\sigma_k(y-\cdot)H(\cdot)(y-\cdot)\right|^q\,dy
\right]^{\frac{1}{q}}\right\|_X\\
&\quad\leq\limsup_{r\rightarrow0^+}\left\|
\left[\fint_{B(\cdot,r)}
\frac{|f(y)-T_{\cdot}^{k-1}(f)(y)-\frac{1}{k!}
\sigma_k(y-\cdot)H(\cdot)(y-\cdot)|^q}{|y-\cdot|^{kq}}
\,dy\right]^{\frac{1}{q}}\right\|_X=0.
\end{align*}
From Lemma \ref{yiyou}, it follows that the sufficiency holds.

Next, we prove the necessity (``only if'' part).
For any $f\in W^{k,X}$ and almost every $x\in\mathbb{R}^n$,
let $H(x):=\nabla(D^{k-1}f)(x)$. Then
$|H|\in X$. Thus, we only need to
show that, for any $f\in W^{k,X}$,
\begin{align*}
\limsup_{r\to 0^+}\left\|\left[\fint_{B(\cdot,r)}
\frac{|f(y)-T_{\cdot}^{k}(f)(y)|^q}{|\cdot-y|
^{kq}}\,dy\right]^{\frac{1}{q}}
\right\|_{X}=0.
\end{align*}
We first show that
this holds for any $g\in C_{\mathrm{c}}^{\infty}$.
Fix $g\in C_{\mathrm{c}}^{\infty}$. For any $r\in(0,1)$,
let
$\rho_r:=\frac{\mathbf{1}_{B(\mathbf{0},r)}}
{|B(\mathbf{0},r)|}$. Then it is easy to verify that
$\{\rho_{r}\}_{r\in(0,1)}$ satisfies \eqref{rou}
and the conditions of Definition \ref{djh}.
From this and Theorem \ref{taile}, we infer that
\begin{align*}
&\lim_{r\to 0^+}\left\|\left[\int_{\mathbb{R}^n}
\frac{|g(y)-T_{\cdot}^{k}(g)(y)|^q}{|\cdot-y|^{(k+1)q}}
\rho_r(\cdot-y)\,dy\right]^{\frac{1}{q}}
\right\|_{X}\\
&\quad=\frac{1}{(1+k)![\sigma (\mathbb{S}^{n-1})]
^{\frac{1}{q}}}\left\|\left[\int_{\mathbb{S}^{n-1}}
\left|\sigma _{k+1}(\xi ) \nabla \left( D^{k}g \right)
(\cdot  ) \xi \right|^q\,d\sigma (\xi )\right]^{\frac{1}{q}}
\right\|_{X}<\infty.
\end{align*}
Then there exist $G\in(0,\infty)$ and $r_0\in(0,\infty)$
such that, for any $r\in(0,r_0)$,
$$\left\|\left[\int_{\mathbb{R}^n}
\frac{|g(y)-T_{\cdot}^{k}(g)(y)|^q}{|\cdot-y|^{(k+1)q}}
\rho_r(\cdot-y)\,dy\right]^{\frac{1}{q}}
\right\|_{X}\leq G.$$
By this, we obtain
\begin{align}\label{ccinfty}
&\limsup_{r\to 0^+}\left\|\left[\fint_{B(\cdot,r)}
\frac{|g(y)-T_{\cdot}^{k}(g)(y)|^q}{|\cdot-y|
^{kq}}\,dy\right]^{\frac{1}{q}}
\right\|_{X}\notag\\
&\quad=\limsup_{r\to 0^+}\left\|
\left[\fint_{B(\cdot,r)}
\frac{|\cdot-y|^q|g(y)-T_{\cdot}^{k}(g)(y)|^q}
{|\cdot-y|^{(k+1)q}}\,dy\right]^{\frac{1}{q}}
\right\|_{X}\notag\\
&\quad\leq\limsup_{r\to 0^+}r\left\|
\left[\fint_{B(\cdot,r)}
\frac{|g(y)-T_{\cdot}^{k}(g)(y)|^q}
{|\cdot-y|^{(k+1)q}}\,dy\right]^{\frac{1}{q}}
\right\|_{X}\leq G\limsup_{r\to 0^+}r=0.
\end{align}
Next, fix $f\in W^{k,X}$. By the density of
$C_{\rm{c}}^{\infty}$ in $\dot{W}^{k,X}$
(see, for instance, \cite[Theorem 2.16(ii)]{hlyyz}),
we conclude that, for any $\varepsilon\in(0,\infty)$,
there exists $g\in C_{\rm{c}}^{\infty}$
such that
\begin{equation}\label{choumi}
\sum_{\left| \beta \right|= k}\frac{k!}{\beta!}\left\|\,
\left|\partial ^{\beta}( f-g )\right|\,
\right\| _X<\varepsilon.
\end{equation}
From \cite[Proposition 2.12]{hlyyz}, we deduce that
\begin{equation}\label{jubu}
f-g\in  W^{k,X}\subset \dot{W}^{k,X}\subset
W^{k,L^1}_{\rm{loc}}.
\end{equation}
Applying the triangle inequality and the
definition of ball Banach function spaces, we find that
\begin{align*}
\text{O}&:=\limsup_{r\rightarrow0^+}
\left\|\left[\fint_{B(\cdot,r)}
\frac{|f(y)-T_{\cdot}^{k}(f)(y)|^q}
{|y-\cdot|^{kq}}\,dy\right]^{\frac{1}{q}}\right\|_X\\
&\ \leq\limsup_{r\rightarrow0^+}\left\|
\left[\fint_{B(\cdot,r)}
\frac{|g(y)-T_{\cdot}^{k}(g)(y)|^q}{|y-\cdot|^{kq}}\,dy
\right]^{\frac{1}{q}}\right\|_X\\
&\quad+\limsup_{r\rightarrow0^+}
\left\|\left[\fint_{B(\cdot,r)}
\frac{|(f-g)(y)-T_{\cdot}^{k-1}(f-g)(y)|^q}
{|y-\cdot|^{kq}}\,dy\right]^{\frac{1}{q}}\right\|_X\\
&\quad+\limsup_{r\rightarrow0^+}\left\|
\left[\fint_{B(\cdot,r)}\frac{|\sum_{|\beta|=k}
\frac{\partial ^{\beta}(f-g)(\cdot)}{\beta!}
(y-\cdot)^{\beta}|^q}{|y-\cdot|^{kq}}\,dy\right]
^{\frac{1}{q}}\right\|_X.
\end{align*}
From this, \eqref{ccinfty}, \eqref{choumi}, \eqref{jubu},
Lemma \ref{zybds1}, and the assumptions that
$\{\mathcal{B}_r\}_{r\in(0,\infty)}$
is uniformly bounded on $X^{\frac{1}{q}}$ and that
$\mathcal{M}$ is bounded on $X$, we deduce that
\begin{align*}
\text{O}&\lesssim\limsup_{r\rightarrow0^+}\left\|
\left[\fint_{B(\cdot,r)}\left|\mathcal{M}
\left(\left|\nabla^k(f-g)\right|\right)(\cdot)
+\mathcal{M}\left(\left|\nabla^k(f-g)\right|\right)(y)
\right|^q\,dy\right]^{\frac{1}{q}}\right\|_X\\
&\quad+\limsup_{r\rightarrow0^+}\sum_{|\beta|=k}\left\|\,
\left|\partial^{\beta}(f-g)\right|\,\right\|_X\\
&\leq\limsup_{r\rightarrow0^+}\left\|\left[
\fint_{B(\cdot,r)}\left|\mathcal{M}\left(
\left|\nabla^k(f-g)\right|\right)(y)\right|^q\,dy\right]
^{\frac{1}{q}}\right\|_X
+\left\|\mathcal{M}\left(\left|\nabla^k(f-g)\right|
\right)\right\|_X+\varepsilon\\
&\lesssim\limsup_{r\rightarrow0^+}\left\|\fint_{B(\cdot,r)}
\left|\mathcal{M}\left(\left|\nabla^k(f-g)\right|\right)
(y)\right|^q\,dy\right\|_{X^{\frac{1}{q}}}^{\frac{1}{q}}
+\left\|\,\left|\nabla^k(f-g)\right|\, \right\|_X
+\varepsilon\\
&\lesssim\left\|\,\left|\nabla^k(f-g)\right|\,
\right\|_X+\varepsilon\lesssim\sum_{|\beta|=k}\frac{k!}{\beta!}\left\|
\,\left|\partial^{\beta}(f-g)\right|\,\right\|_X
+\varepsilon\lesssim\varepsilon.
\end{align*}
By the arbitrariness of $\varepsilon\in(0,\infty)$,
we obtain
$$\lim_{r\rightarrow0^+}\left\|
\left[\fint_{B(\cdot,r)}\frac{|f(y)-T_{\cdot}^{k}(f)(y)|^q}
{|y-\cdot|^{kq}}\,dy\right]^{\frac{1}{q}}\right\|_X=0.$$
This completes the proof of Theorem \ref{tailek}.
\end{proof}
\section{Proof of Theorem \ref{the-4}}\label{s7}
In this section, we prove Theorem \ref{the-4}, which
gives a negative answer to the question in
\cite[Remark 3.19(ii)]{jyy25}, which asks whether
$f\in W^{k,X}$ if and only if $f\in W^{k-1,X}$ and
\begin{align}\label{equ-local-van}
\lim_{r\to 0+}\left\|\frac{1}{r^k}
\left[\fint_{B(\cdot,r)}\left|f(y)-P_{B(\cdot,r)}
^{(k)}(f)(y)\right|^q\,dy\right]
^{\frac{1}{q}}\right\|_X=0.
\end{align}
To prove Theorem \ref{the-4}, we take $q:=2$, $k=1$, and $X:=L^2(\mathbb{R})$.
Inspired by the estimate for finite
differences in \cite[Theorem 2.1]{gt},
we first establish
the following estimate for minimizing polynomials.
\begin{lemma}\label{zhengjiaotouying}
There exists a
positive constant $C$ such that, for any $f\in L^2(\mathbb{R})$
and $r\in(0,\infty)$,
\begin{align}\label{equ-four-f}
\left\|\left[\fint_{B(\cdot,r)}
\left|f(y)-P_{B(\cdot,r)}^{(1)}(f)(y)\right|^2\,dy\right]
^{\frac{1}{2}}
\right\|_{L^2(\mathbb{R})}\leq C\left[\int_{\mathbb{R}}
\left|\mathcal{F}f(\xi )\right|^2
\min\left\{1,|r\xi|^4\right\}\,d\xi \right]^{\frac{1}{2}}.
\end{align}
\end{lemma}
To prove this lemma, we need the following
inequality, which is precisely
\cite[(2.12)]{jtyyz}.
\begin{lemma}\label{jxyzx}
Let $q\in[1,\infty)$ and $k\in \mathbb{Z}_+$.
Then,
for any complex-valued function
$f\in L_{\mathrm{loc}}^1$ and any ball $B\subset\mathbb{R}^n$,
$$\left[\fint_{B}\left|f(x)-P_{B}^{(k)}(f)(x)\right|^q\,dx\right]
^{\frac{1}{q}}\sim\inf_{P\in\mathcal{P}_k}
\left[\fint_{B}\left|f(x)-P(x)\right|^q\,dx\right]
^{\frac{1}{q}},$$
where the positive equivalence constants are independent of
$B$ and $f$. The same estimate
holds with $B$ replaced by any cube $Q$.
\end{lemma}
\begin{proof}[Proof of Lemma \ref{zhengjiaotouying}]
Fix
$f\in L^1(\mathbb{R})\cap L^2(\mathbb{R})$, $x\in\mathbb{R}$, and $r\in(0,\infty)$.
From a change of variables, we infer that
\begin{align}\label{equ-fint}
\fint_{B(x,r)}
\left|f(y)-P_{B(x,r)}^{(1)}(f)(y)\right|^2\,dy
&=\fint_{B(\mathbf{0},1)}
\left|f(x+rt)-P_{B(x,r)}^{(1)}(f)(x+rt)\right|^2dt.
\end{align}
Let $\{\phi_{\beta}:\beta\in\{0,1\}\}$
be the \emph{orthogonal polynomials} obtained from
$\{x^\beta:\beta\in\{0,1\}\}$ on $B(\mathbf{0},1)$ by the
Gram--Schmidt method with respect
to the weight $\frac{1}{|B(\mathbf{0},1)|}$.
That is, $\{\phi_\beta:\beta\in\{0,1\}\}\subset\mathcal{P}_1$
and, for any $\beta_1,\beta_2\in\{0,1\}$,
\begin{align*}
\left(\phi_{\beta_1},\phi_{\beta_2}\right)_{B(\mathbf{0},1)}
:=\fint_{B(\mathbf{0},1)}
{\phi_{\beta_1}(x)\overline{\phi_{\beta_2}(x)}\,dx}
=\begin{cases}
	1,&		\beta _1=\beta _2,\\
	0,&		\beta _1\ne \beta _2.\\
\end{cases}
\end{align*}
It is well known that, for any $f\in L_{\mathrm{loc}}^1$
and $y\in\mathbb{R}$,
\begin{align}\label{equ-pol-f}
P_{B(x,r)}^{(1)}(f)(y)
=\sum_{\beta=0}^{1}
\fint_{B(x,r)}
\overline{\phi_{\beta}\left(\frac{z-x}{r}\right)}
f(z)\,dz\phi_{\beta}\left(\frac{y-x}{r}\right);
\end{align}
see, for instance, \cite[p.\,55,\ line 16]{l}.
By this, a change of variables, and the fact that
$f\in L^1(\mathbb{R})$
and $\phi_\beta\in L^1[-1,1]$ for any
$\beta\in\{0,1\}$,
we conclude that, for any $t\in[-1,1]$,
\begin{align}\label{equ-pol-f-L1}
\left\|P_{B(\cdot,r)}^{(1)}(f)(\cdot+rt)
\right\|_{L^1(\mathbb{R})}
&\leq \sum_{\beta=0}^1\left\|\fint_{B(\mathbf{0},1)}
\overline{\phi_{\beta}(h)}f(\cdot+rh)\,dh\phi_{\beta}(t)
\right\|_{L^1(\mathbb{R})}\notag\\
&\leq\sum_{\beta=0}^1\fint_{B(\mathbf{0},1)}|
\phi_{\beta}(h)|\left\|f(\cdot+rh)
\right\|_{L^1(\mathbb{R})}\,dh|\phi_{\beta}(t)|\notag\\
&=\left\|f\right\|_{L^1(\mathbb{R})}\sum_{\beta=0}^1\fint_{B(\mathbf{0},1 )}|
\phi_{\beta}(h)|\,dh|\phi_{\beta}(t)|<\infty,
\end{align}
and hence $P_{B(\cdot,r)}^{(1)}(f)(\cdot+rt)\in
L^1(\mathbb{R})$. For any $\xi\in\mathbb{R}$ and $t\in[-1,1]$,
define $g_\xi(t):=e^{ir\xi t}$.
By this, \eqref{equ-pol-f}, Fubini's theorem,
and a change of variables,
we conclude that, for any $\xi\in\mathbb{R}$ and $t\in[-1,1]$,
\begin{align}\label{equ-plancherel}
\mathcal{F}\left(P_{B(\cdot,r)}^{(1)}(f)(\cdot+rt)\right)(\xi)&=
\frac{1}{\sqrt{2\pi }}\int_{\mathbb{R}}e^{-ix\xi}P_{B(x,r)}
^{(1)}(f)(x+rt)\,dx\notag\\
&=\sum_{\beta=0}^1\frac{1}{\sqrt{2\pi }}
\int_{\mathbb{R}}e^{-ix\xi}\fint_{B( x,r )}
\overline{\phi _{\beta}\left(\frac{z-x}{r}\right)}f(z)\,dz
\phi _{\beta}\left(\frac{x+rt-x}{r}\right)\,dx\notag\\
&=\sum_{\beta=0}^{1}
\fint_{B(\mathbf{0},1 )}
\overline{\phi _{\beta}(h)}
\frac{1}{\sqrt{2\pi }}
\int_{\mathbb{R}}e^{-ix\xi}f(x+rh)\,dx
\,dh\phi _{\beta}(t)\notag\\
&=\sum_{\beta=0}^1
\fint_{B(\mathbf{0},1)}
\overline{\phi_{\beta}(h)}e^{irh\xi}
\frac{1}{\sqrt{2\pi}}
\int_{\mathbb{R}}e^{-i(x+rh)\xi}f(x+rh)\,dx
\,dh\phi_{\beta}(t)\notag\\
&=\mathcal{F}f(\xi)
\sum_{\beta=0}^1
\fint_{B(\mathbf{0},1)}\overline{\phi_{\beta}(h)}e^{irh\xi}\,dh
\phi_{\beta}(t)=\mathcal{F}f(\xi)P_{B(\mathbf{0},1)}
^{(1)}(g_{\xi})(t).
\end{align}
Applying an argument similar to that used in the estimation of
\eqref{equ-pol-f-L1} with $L^1$ replaced by $L^2$, and applying
Minkowski's inequality,
we conclude that, for any $t\in[-1,1]$,
$P_{B(\cdot,r)}^{(1)}(f)(\cdot+rt)\in L^2$.
From this, \eqref{equ-fint},
\eqref{equ-plancherel}, Tonelli's theorem, and Plancherel's theorem, we infer that
\begin{align}\label{equ-f-Pol}
&\left\|
\left[\fint_{B(\cdot,r)}
\left|f(y)-P_{B(\cdot,r)}^{(1)}(f)(y)\right|^2\,dy\right]
^{\frac{1}{2}}
\right\|_{L^2(\mathbb{R})}^2\notag\\
&\quad=\fint_{B(\mathbf{0},1)}
\int_{\mathbb{R}}
\left|f(x+rt)-P_{B(x,r)}^{(1)}(f)(x+rt)\right|^2\,dx\,dt\notag\\
&\quad=\fint_{B(\mathbf{0},1)}
\int_{\mathbb{R}}
\left|\mathcal{F}\left(f(\cdot+rt)\right)(\xi)-\mathcal{F}
\left(P_{B(\cdot,r)}^{(1)}(f)(\cdot+rt)\right)(\xi)
\right|^2\,d\xi\,dt\notag\\
&\quad=\fint_{B(\mathbf{0},1)}
\int_{\mathbb{R}}
\left|g_{\xi}(t)-P_{B(\mathbf{0},1)}
^{(1)}(g_{\xi})(t)\right|^2|\mathcal{F}{f}(\xi)|^2\,d\xi\,
dt\notag\\
&\quad=\int_{\mathbb{R}}
|\mathcal{F}{f}(\xi)|^2\fint_{B(\mathbf{0},1)}
\left|g_{\xi}(t)-P_{B(\mathbf{0},1)}^{(1)}
(g_{\xi})(t)\right|^2\,dt\,d\xi .
\end{align}
By Lemma \ref{jxyzx}, we obtain
\begin{align}\label{equ-gxi-rxi}
\fint_{B(\mathbf{0},1)}
\left|g_{\xi}(t)-P_{B(\mathbf{0},1)}
^{(1)}(g_{\xi})(t)\right|^2\,dt&\lesssim
\fint_{B(\mathbf{0},1)}
\left|e^{ir\xi t}-1-ir\xi t\right|^2\,dt\notag\\
&=\fint_{B(\mathbf{0},1)}
\left|\int_{0}^{r\xi t}(r\xi t-s)e^{is}\,ds\right|^2\,dt\notag\\
&\le\fint_{B(\mathbf{0},1)}
\left[\int_{0}^{|r\xi t|}
(|r\xi t|-s)\,ds\right]^2\,dt\sim|r\xi|^4,
\end{align}
with the implicit positive constants independent of $r,\xi$, and $t$.
On the other hand, from Lemma \ref{jxyzx}, we deduce that
\begin{align*}
\fint_{B(\mathbf{0},1)}
\left|g_{\xi}(t)-P_{B(\mathbf{0},1)}^{(1)}(g_{\xi})(t)
\right|^2\,dt\lesssim
\fint_{B(\mathbf{0},1)}
\left|e^{ir\xi t}\right|^2\,dt\sim 1,
\end{align*}
where the implicit positive constant is independent
of $r,\xi$, and $t$.
This, together with \eqref{equ-f-Pol}
and \eqref{equ-gxi-rxi},
further implies that
\begin{align*}
\left\|\left[\fint_{B(\cdot,r)}
\left|f(y)-P_{B(\cdot,r)}^{(1)}(f)(y)\right|^2\,dy\right]
^{\frac{1}{2}}
\right\|_{L^2(\mathbb{R})}
\lesssim\left[\int_{\mathbb{R}}|\mathcal{F}f(\xi )|^2
\min\left\{1,|r\xi|^4\right\}\,d\xi \right]^{\frac{1}{2}},
\end{align*}
where the implicit positive constant is
independent of $r$ and $f$.
Thus, we show that
\eqref{equ-four-f} holds for any $f\in L^1(\mathbb{R})\cap L^2(\mathbb{R})$.
On the other hand, by Lemma \ref{jxyzx} and Fubini's theorem,
we conclude that, for any
$f\in L^2(\mathbb{R})$,
\begin{align}\label{equ-L2bound}
\left\|\left[\fint_{B(\cdot,r)}
\left|f(y)-P_{B(\cdot,r)}^{(1)}(f)(y)\right|^2\,dy\right]
^{\frac{1}{2}}
\right\|_{L^2(\mathbb{R})}\lesssim\left\|
\left[\fint_{B(\cdot,r)}
|f(y)|^2\,dy\right]^{\frac{1}{2}}
\right\|_{L^2(\mathbb{R})}=
\left\|f\right\|_{L^2(\mathbb{R})}.
\end{align}
Fix $f\in L^2(\mathbb{R})$.
Since $L^1(\mathbb{R})\cap L^2(\mathbb{R})$ is dense
in $L^2(\mathbb{R})$,
it follows that there exists $\{f_j\}_{j\in\mathbb{N}}$
in $L^1(\mathbb{R})\cap L^2(\mathbb{R})$
such that $f_j\to f$ in $L^2(\mathbb{R})$
as $j\to\infty$.
Therefore, using \eqref{equ-L2bound}, we obtain
\begin{align}\label{equ-f-polf-l2}
\lim_{j\to\infty}\left\|\left[\fint_{B(\cdot,r)}
\left|(f-f_j)(y)-P_{B(\cdot,r)}^{(1)}(f-f_j)(y)\right|^2
\,dy\right]^{\frac{1}{2}}
\right\|_{L^2(\mathbb{R})}\lesssim
\lim_{j\to\infty}\left\|f-f_j\right\|_{L^2(\mathbb{R})}=0.
\end{align}
On the other hand, from Plancherel's theorem,
we deduce that
\begin{align*}
\lim_{j\to\infty}
\left[\int_{\mathbb{R}}
|(\mathcal{F}{f}-\mathcal{F}{f_j})(\xi )|^2
\min\left\{1,|r\xi|^4\right\}\,d\xi \right]^{\frac{1}{2}}
&\leq\lim_{j\to\infty}\left[\int_{\mathbb{R}}
|\mathcal{F}(f-f_j)(\xi )|^2\,d\xi \right]^{\frac{1}{2}}\\
&=\lim_{j\to\infty}\left\|f-f_j\right\|_{L^2(\mathbb{R})}=0.
\end{align*}
This, together with \eqref{equ-f-polf-l2}, further implies that
\begin{align*}
\left\|
\left[\fint_{B(\cdot,r)}
\left|f(y)-P_{B(\cdot,r)}^{(1)}(f)(y)\right|^2\,dy\right]
^{\frac{1}{2}}
\right\|_{L^2(\mathbb{R})}
&=\lim_{j\to\infty}
\left\|
\left[\fint_{B(\cdot,r)}
\left|f_j(y)-P_{B(\cdot,r)}^{(1)}(f_j)(y)\right|^2\,dy\right]
^{\frac{1}{2}}
\right\|_{L^2(\mathbb{R})}\\
&\lesssim
\lim_{j\to\infty}\left[\int_{\mathbb{R}}|\mathcal{F}f_j(\xi )|^2
\min\left\{1,|r\xi|^4\right\}\,d\xi \right]^{\frac{1}{2}}\\
&=\left[\int_{\mathbb{R}}|\mathcal{F}f(\xi )|^2
\min\left\{1,|r\xi|^4\right\}\,d\xi \right]^{\frac{1}{2}},
\end{align*}
which shows \eqref{equ-four-f}. This completes the proof
of Lemma \ref{zhengjiaotouying}.
\end{proof}
Next, with the aid of Lemma \ref{zhengjiaotouying},
we can show Theorem \ref{the-4}.
\begin{proof}[Proof of Theorem \ref{the-4}]
For any $\xi\in\mathbb{R}$, let
\begin{align}\label{equ-def-f}
\mathcal{F}f( \xi ) :=\sum_{j=2}^{\infty}
\frac{1}{2^{j}\sqrt{j}}
\left[\mathbf{1}_{[ 2^j,2^j+1)}(\xi)+
\mathbf{1}_{( -2^j-1,-2^j]}(\xi) \right],
\end{align}
and hence
\begin{align}\label{equ-Ff-L2}
\|\mathcal{F}{f}\|_{L^2(\mathbb{R})}
=\left(2\sum_{j=2}^{\infty}\frac{1}{2^{2j}j}\right)^{\frac{1}{2}}<\infty.
\end{align}
From \cite[p.\,22]{gtm250} and \eqref{equ-def-f}, we deduce that
\begin{align*}
\|\nabla f\|_{L^2(\mathbb{R})}^2&\sim
||i(\cdot)\mathcal{F}(f)(\cdot)||_{L^2(\mathbb{R})}^2
=\int_{\mathbb{R}}|\xi|^2|\mathcal{F}{f}(\xi)|^2\,d\xi\\
&=\sum_{j=2}^{\infty}
\int_{2^j\le|\xi|< 2^{j}+1}
|\xi|^2\left|\mathcal{F}f(\xi)\right|^2\,d\xi
\geq2\sum_{j=2}^{\infty}\frac{2^{2j}}{2^{2j}j}=\infty,
\end{align*}
which, together with \eqref{equ-Ff-L2},
further implies that
$f\in L^2(\mathbb{R})\setminus W^{1,2}(\mathbb{R})$.
We now prove
$$\lim_{r\to 0^+}
\left\|\frac{1}{r}
\left[\fint_{B(\cdot,r)}
\left|f(y)-P_{B(\cdot,r)}
^{(1)}(f)(y)\right|^2\,dy\right]^{\frac{1}{2}}
\right\|_{L^2(\mathbb{R})}=0.$$
Fix $r\in(0,\infty)$.
Using Lemma \ref{zhengjiaotouying}
and the pairwise disjointness of
$\{[2^{j},2^{j}+1),(-2^{j}-1,-2^{j}]\}_{j=2}^\infty$, we obtain
\begin{align}\label{equ-f-pol-four}
&\left\|
\left[\fint_{B(\cdot,r)}\left|
f(y)-P_{B(\cdot,r)}^{(1)}(f)(y)\right|^2\,dy\right]^{\frac{1}{2}}
\right\|_{L^2(\mathbb{R})}\notag\\
&\quad\lesssim
\left[\int_{\mathbb{R}}|\mathcal{F}{f}(\xi )|^2
\min\left\{1,|r\xi|^4\right\}\,d\xi \right]^{\frac{1}{2}}\notag\\
&\quad=\left[\int_{\mathbb{R}}
\sum_{j=2}^{\infty}\frac{1}{4^{j}j}
\left[\mathbf{1}_{[ 2^j,2^j+1 )}(\xi)
+\mathbf{1}_{( -2^j-1,-2^j]}(\xi)\right]
\min\left\{1,|r\xi|^4\right\}\,d\xi\right]^{\frac{1}{2}}.
\end{align}
Fix $r\in(0,\frac{1}{8})$, and let $k\in\mathbb{N}$ be
such that $2^{-k}<r\leq2^{-k+1}$. From this, we infer that
\begin{align}\label{equ-I-II-r}
&\int_{\mathbb{R}}\sum_{j=2}^{\infty}\frac{1}{4^{j}j}
\left[\mathbf{1}_{[ 2^j,2^j+1 )}(\xi)+
\mathbf{1}_{( -2^j-1,-2^j]}(\xi) \right]
\min\left\{1,|r\xi|^4\right\}\,d\xi\notag\\
&\quad=\sum_{j=2}^{\infty}\frac{2}{4^{j}j}
\int_{2^j}^{2^j+1}\min\left\{1,|r\xi|^4\right\}\,d\xi
\leq\sum_{j=2}^{k-2}\frac{2}{4^{j}j}
\int_{2^j}^{2^j+1}|r\xi|^4\,d\xi
+\sum_{j=k-1}^{\infty}\frac{2}{4^{j}j}
\int_{2^j}^{2^j+1}\,d\xi\notag\\
&\quad=:\mathrm{I}(r)+\mathrm{II}(r).
\end{align}
We first estimate $\mathrm{I}(r)$.
Then we have
\begin{align*}
\mathrm{I}(r)
\sim\sum_{j=2}^{k-2}
\frac{4^{2j}r^4}{4^{j}j}=\sum_{j=2}^{k-2}
\frac{4^{j}r^4}{j}
\lesssim r^4\sum_{j=2}^{k-2}\int_{j}^{j+1}\frac{4^x}{x}\,dx=r^4
\int_{2}^{k-1}
\frac{4^x}{x}\,dx\le r^4\int_{2}^{-\frac{\log r}{\log 2}}
\frac{4^x}{x}\,dx.
\end{align*}
This further implies that
\begin{align}\label{equ-Ir}
\lim_{r\to 0^+}r^{-2}\mathrm{I}(r)
=\lim_{r\to 0^+}r^{2}
\int_{2}^{-\frac{\log r}{\log 2}}
\frac{4^x}{x}\,dx=0.
\end{align}
On the other hand, we have
\begin{align*}
\mathrm{II}(r)
\le
\sum_{j=k-1}^{\infty}
\int_{j-1}^j\frac{1}{4^{x}x}\,dx
=\int_{k-2}^\infty\frac{1}{4^{x}x}\,dx
\le
\int_{-\frac{\log r}{\log 2}-2}^\infty\frac{1}{4^{x}x}\,dx,
\end{align*}
which further implies that
\begin{align*}
\lim_{r\to 0^+}
r^{-2}\mathrm{II}(r)
\le\lim_{r\to 0^+}
r^{-2}\int_{-\frac{\log r}{\log 2}-2}^\infty\frac{1}{4^{x}x}\,dx=0.
\end{align*}
This, together with \eqref{equ-f-pol-four}, \eqref{equ-I-II-r},
and \eqref{equ-Ir},
further implies that
\begin{align*}
\lim_{r\to 0^+}
\left\|\frac{1}{r}
\left[\fint_{B(\cdot,r)}
|f(y)-P_{B(\cdot,r)}^{(1)}(f)(y)|^2\,dy\right]^{\frac{1}{2}}
\right\|_{L^2(\mathbb{R})}^2
\lesssim
\lim_{r\to 0^+}
r^{-2}\left[\mathrm{I}(r)+\mathrm{II}(r)\right]=0.
\end{align*}
This completes the proof of Theorem \ref{the-4}.
\end{proof}
\section{Proofs of Theorems \ref{pro-laplace} and \ref{pro-weaknorm}}\label{s8}
This section is divided into two subsections.
In Subsection \ref{sub-laplace}, we disprove the
proposed BBM-type formula involving ball averages in \cite[Remark 4.10]{jyy25}
(see Theorem \ref{pro-laplace}). In Subsection \ref{sub-weak-laplace},
we show that its BSVY-type analogue also fails
(see Theorem \ref{pro-weaknorm}). Both results are established
by using the singular radial kernels constructed in Section \ref{s4}.
\subsection{Proof of Theorem \ref{pro-laplace}}\label{sub-laplace}
In this subsection, we prove Theorem \ref{pro-laplace}.
To this end, we first recall the question proposed in
\cite[Remark 4.10]{jyy25}.
Let $X$ be a ball Banach function space satisfying some mild
assumptions (see \cite[Theorem 4.4]{jyy25}). Let $q\in[1,\infty)$ and
$k\in\mathbb{N}$.
The question is whether, for any $f\in C_{\mathrm{c}}^{2k}$,
\begin{equation}\label{equ-bbm-sim}
\lim_{r\to 0}
\left\|\left[\int_{\mathbb{R} ^n}
{\frac{|f(y)-B_{k,r}(f)
(y) |^q}{|\cdot-y|^{2kq}}}\rho _r(\cdot-y)
\,dy\right]^{\frac{1}{q}}\right\|
_X\sim \left\|\Delta^{k}f\right\|_X,
\end{equation}
where $\{\rho_r\}_{r\in(0,\infty)}$ satisfies the conditions of Definition
\ref{djh}, $B_{k,r}$ is the higher-order average operator
(see \cite[(3.3)]{dgyy15} for the precise definition),
and the positive equivalence constants are
independent of $f$.

Theorem \ref{pro-laplace} gives a negative answer to this question.
Specifically, let $q\in(1,\infty)$, $X=L^q$, and $k=1$.
Then, by the definition of the higher-order average
operator,
we find that,
for any $f\in L_{\mathrm{loc}}^1$, $r\in(0,\infty)$, and
$y\in\mathbb{R}^n$, $B_{1,r}(f)(y)=f_{B(y,r)}$.
We now prove Theorem \ref{pro-laplace}.

\begin{proof}[Proof of Theorem \ref{pro-laplace}]
From \eqref{equ-rho-1} and \eqref{equ-rho-0},
we deduce that $\{\rho_r\}_{r\in(0,\infty)}$
satisfies the conditions of Definition \ref{djh}.
The sufficiency (``if'' part) is easily deduced
from the mean-value property of
harmonic functions
(see, for instance, \cite[Theorem 2]{e10}), and
hence we only need to show the necessity (``only if'' part).
Fix $f\in C_{\mathrm{c}}^{3}$ satisfying \eqref{equ-bbm}.
Since $f\in C_{\mathrm{c}}^3$, there exists
$M\in(0,\infty)$ such that, for any $x\in\mathbb{R}^n$,
\begin{equation}\label{equ-3d}
\sum_{\{\gamma\in\mathbb{Z}_+^n:|\gamma|=3\}}
\left|\frac{\partial^{\gamma}f(x)}{\gamma !}\right|<M.
\end{equation}
Next, we show that $\Delta f\equiv0$.
Otherwise, there exists
$x_0\in\mathbb{R}^n$ such that $|\Delta f(x_0)|>0$, and hence
there exists $\eta\in(0,\infty)$ such that,
for any $x\in B(x_0,\eta)$,
$|\Delta f(x)|>\frac{|\Delta f(x_0)|}{2}$,
which further implies that
\begin{equation}\label{equ-inf}
\alpha:=\inf_{x\in B(x_0,\eta)}
\frac{|\Delta f(x)|}{6M(n+2)}>0.
\end{equation}
On the other hand, by the Taylor remainder theorem,
we conclude that,
for any $x,y\in\mathbb{R}^n$, there exists $\vartheta\in(0,1)$
such that
\begin{equation}\label{equ-taylor}
f(y)
=\sum_{\{\gamma\in\mathbb{Z}_+^n,|\gamma|\le 2\}}
\frac{\partial^\gamma f(x)}{\gamma !}(y-x)^{\gamma}
+\mathcal{R}^{(2)}(x,y-x),
\end{equation}
where
\begin{equation*}
\mathcal{R}^{(2)}(x,y-x):=
\sum_{\{\gamma\in\mathbb{Z}_+^n:|\gamma|=3\}}
\frac{\partial^\gamma f(\vartheta x+(1-\vartheta)y)}
{\gamma !}(y-x)^{\gamma}.
\end{equation*}
From this and \eqref{equ-3d}, we deduce that,
for any $r\in(0,\infty)$ and $y\in B(x,r)$,
\begin{equation}\label{equ-tayre}
\left|\mathcal{R}^{(2)}(x,y-x)\right|\le Mr^3.
\end{equation}
Using \eqref{teshu}, \eqref{equ-bbm}, and
\eqref{equ-inf}, we conclude that
there exists $r_0\in (0,\frac{\alpha}{2})$ such that
\begin{equation}\label{equ-leinfty}
\frac{\theta}{r_0^{\theta}\sigma(\mathbb{S}^{n-1})}\int_{\mathbb{R}^n}
{\int_{B(x,r_0)}{\frac{|f(y)-f_{B(y,r_0)}|^q}
{|x-y|^{2q+n-\theta}}\,dy}}\,dx<\infty.
\end{equation}
On the other hand, averaging both sides of \eqref{equ-taylor}
over $y\in B(x,r_0)$, we find that,
for any $x\in B(x_0,\eta)$,
\begin{equation}\label{equ-mean}
f(x)-f_{B(x,r_0)}
=-\frac{1}{2(n+2)}\Delta f(x) r_0^2
-\fint_{B(x,r_0)}\mathcal{R}^{(2)}(x,y-x)\,dy.
\end{equation}
Using a change of variables and Tonelli's theorem,
we conclude that
\begin{align}\label{equ-var}
&\int_{\mathbb{R} ^n}{\int_{B(x,r_0)}
{\frac{|f(y)-f_{B(y,r_0)}|^q}{|x-y|^{2q+n-\theta}}\,dy}}\,
dx\notag\\
&\quad
=\int_{\mathbb{R} ^n}{\int_{B(\mathbf{0},r_0)}
{\frac{|f(x-h)-f_{B(x-h,r_0)}|^q}{|h|^{2q+n-\theta}}\,dh}}\,dx\notag\\
&\quad
\sim \sum_{j=0}^{\infty}2^{j(n+2q-\theta)}
r_0^{-(n+2q-\theta)}
\int_{\mathbb{R}^n}
\int_{B(\mathbf{0},2^{-j}r_0)\setminus
B(\mathbf{0},2^{-j-1}r_0)}
|f(x-h)-f_{B(x-h,r_0)}|^q\,dh\,dx\notag \\
&\quad
=\sum_{j=0}^{\infty}2^{j(n+2q-\theta)}r_0^{-(n+2q-\theta)}
\int_{\mathbb{R}^n}
\int_{B(\mathbf{0},2^{-j}r_0)\setminus
B(\mathbf{0},2^{-j-1}r_0)}
|f(x)-f_{B(x,r_0)}|^q\,dh\,dx\notag\\
&\quad
\sim\sum_{j=0}^{\infty}2^{j(2q-\theta)}r_0^{-(2q-\theta)}
\int_{\mathbb{R}^n}
|f(x)-f_{B(x,r_0)}|^q\,dx\notag\\
&\quad
\ge\sum_{j=0}^{\infty}2^{j(2q-\theta)}r_0^{-(2q-\theta)}
\int_{B(x_0,\eta)}
|f(x)-f_{B(x,r_0)}|^q\,dx.
\end{align}
By \eqref{equ-inf},
\eqref{equ-tayre}, \eqref{equ-mean}, the fact
that $r_0\in (0,\frac{\alpha}{2})$,
and Lemma \ref{fmsbds} with $\delta=\frac{1}{2}$
therein,
we conclude that
\begin{align*}
&\int_{B(x_0,\eta)}
|f(x)-f_{B(x,r_0)}|^q\,dx\\
&\quad\geq
\int_{B(x_0,\eta)}
\left[
\frac{1}{(1+\frac{1}{2})^q}\frac{1}{2^q(n+2)^q}
|\Delta f(x) r_0^2|^q-\frac{1}{(\frac{1}{2})^q}
\left|\fint_{B(x,r_0)}\mathcal{R}^{(2)}(x,y-x)\,dy\right|^q
\right]\,dx\\
&\quad\ge
\int_{B(x_0,\eta)}
\left[\frac{1}{3^q(n+2)^q}
|\Delta f(x) r_0^2|^q-2^q
M^q r_0^{3q}\right]\,dx\\
&\quad\ge\int_{B(x_0,\eta)}2^qr_0^{2q}M^q
\left(\alpha^q-r_0^q\right)\,dx\sim 2^q\eta^nr_0^{2q}M^q
\left(\alpha^q-r_0^q\right)>0.
\end{align*}
This, together with \eqref{equ-var} and $\theta<2q$,
further implies that
\begin{align*}
\int_{\mathbb{R} ^n}{\int_{B(x,r_0)}
{\frac{|f(y)-f_{B(y,r_0)}|^q}{|x-y|^{2q+n-\theta}}\,dy}}\,dx
\gtrsim
2^q\eta^nr_0^{\theta}M^q
\left(\alpha^q-r_0^q\right)\sum_{j=0}^{\infty}2^{j(2q-\theta)}
=\infty,
\end{align*}
which contradicts \eqref{equ-leinfty}, and hence $\Delta f\equiv0$. This
completes the proof of Theorem \ref{pro-laplace}.
\end{proof}

\subsection{Proof of Theorem \ref{pro-weaknorm}}\label{sub-weak-laplace}
In this subsection, we prove Theorem \ref{pro-weaknorm},
which shows that the BSVY-type analogue of \eqref{equ-bbm-sim}
also fails in general. For convenience, we first
introduce the corresponding level sets.
For any $r\in(0,\infty)$,
$f\in L_{\rm{loc}}^1$, and $\lambda\in(0,\infty)$, let
\begin{align}\label{Eru}
E_{\lambda,r}(f):=\left\{
(x,y)\in\mathbb{R}^n\times\mathbb{R}^n
:x\ne y,\frac{|f(y)-f_{B(y,r)}|}{|x-y|^{2}}>\lambda
\right\}.
\end{align}
We now prove Theorem \ref{pro-weaknorm}.
\begin{proof}[Proof of Theorem \ref{pro-weaknorm}]
The sufficiency (``if'' part) is obvious,
and hence we only need to show the necessity (``only if'' part).
Fix $f\in C_{\mathrm{c}}^2$.
If \eqref{equ-weakLq} holds, then there exists
$r\in(0,\infty)$ such that
\begin{align}\label{equ-weak-r}
\sup_{\lambda\in(0,\infty)}
\lambda\left[\nu_r(E_{\lambda,r}(f))\right]^{\frac{1}{q}}<\infty.
\end{align}
If $f\not\equiv0$,
then, applying an argument similar to that used in
the estimation of \eqref{equ-f-fav},
we conclude that
there exists $x_0\in\mathbb{R}^n$ such that
$|f(x_0)-f_{B(x_0,r)}|=:\alpha>0$.
By the assumption that $f\in C_{\mathrm{c}}^2$ and
the Lebesgue dominated convergence theorem,
we conclude that
$f(\cdot)-f_{B(\cdot,r)}$ is continuous on $\mathbb{R}^n$.
Therefore, there exists $\eta\in(0,\infty)$
such that, for any $y\in B(x_0,\eta)$,
$|f(y)-f_{B(y,r)}|>\frac{\alpha}{2}$.
From this, we deduce that, for
any $\lambda\in(\frac{\alpha}{2r^2},\infty)$,
$y\in B(x_0,\eta)$, and
$x\in B(y,(\frac{\alpha}{2\lambda})^{\frac{1}{2}})$,
\begin{align*}
\frac{|f(y)-f_{B(y,r)}|}{|x-y|^2}>
\frac{2\alpha \lambda}{2\alpha}=\lambda,
\end{align*}
and hence
$$\left\{(x,y)\in\mathbb{R}^{2n}:
y\in B(x_0,\eta), x\in
B\left(y,\left(\frac{\alpha}{2\lambda}\right)^{\frac{1}{2}}
\right)\right\}
\subset E_{\lambda,r}(f).$$
Using this, we conclude that
\begin{align*}
\nu_r(E_{\lambda,r}(f))&
\ge \int_{B(x_0,\eta )}
\int_{B(y,(\frac{\alpha}{2\lambda})^{\frac{1}{2}})}
\frac{\theta \mathbf{1}_{(0,r)}(|x-y|)}
{r^{\theta}\sigma(\mathbb{S}^{n-1})|x-y|^{n-\theta}}\,dx\,dy\\
&=\int_{B(x_0,\eta)}
\int_{B(\mathbf{0},(\frac{\alpha}{2\lambda})^{\frac{1}{2}})}
\frac{\theta}
{r^{\theta}\sigma(\mathbb{S}^{n-1})|z|^{n-\theta}}\,dz\,dy
\sim
\lambda^{-\frac{\theta}{2}},
\end{align*}
with the implicit positive constant
independent of $\lambda$.
This, together with the arbitrariness of
$\lambda\in(\frac{\alpha}{2r^2},\infty)$
and the assumption that
$\theta\in(0,2q)$, further implies that
\begin{align*}
\sup_{\lambda\in(0,\infty)}
\lambda
\left[\nu_{r}(E_{\lambda,r}(f))\right]^{\frac{1}{q}}
\gtrsim
\sup_{\lambda\in(\frac{\alpha}{2r^2},\infty)}
\lambda^{1-\frac{\theta}{2q}}=\infty,
\end{align*}
which contradicts \eqref{equ-weak-r}.
Therefore, $f\equiv 0$,
which then completes the proof of Theorem \ref{pro-weaknorm}.
\end{proof}
\section{Proof of Theorem \ref{q-wug}}\label{s9}
In this final section, we consider the dependence of $JN^{\rm{con}}_{(p,q,k)_\alpha}$
on the integrability exponent $q$.
To this end, we first recall the
definition of $JN^{\rm{con}}_{(p,q,k)_\alpha}$ (see \cite[Definition 1.5]{jtyyz}
or \cite[Definition 5.17]{jyy25}).
\begin{definition}\label{def-john}
Let $p,q\in[1,\infty)$,
$k\in\mathbb{Z}_+$, and $\alpha
\in\mathbb{R}$. The
\textit{special John--Nirenberg--Campanato space
via congruent cubes}, denoted by
$JN^{\rm{con}}_{(p,q,k)_\alpha}$, is defined as
the set of all $f\in L^q_{\rm{loc}}$ such
that
$$\|f\|_{JN^{\rm{con}}_{(p,q,k)_\alpha}}
:=\sup_{l\in(0,\infty)}[f]_{(p,q,k)_\alpha}^{(l)}<\infty,$$
where, for any $l\in(0,\infty)$,
$$[f]_{(p,q,k)_\alpha}^{(l)}:=\sup_{\{Q_j\}_j\in\Pi_l}
\left\{\sum_{j}\left|Q_j\right|^{1-p\alpha}\left[
\fint_{Q_j}\left|f(x)-P_{Q_j}^{(k)}(f)(x)\right|^q\,dx
\right]^\frac{p}{q}\right\}^{\frac{1}{p}}.$$
Here and in what follows, the \emph{notation $\Pi_l$} denotes the collection of all
families of pairwise disjoint cubes in $\mathbb{R}^n$ with edge length $\ell$.
\end{definition}
To prove Theorem \ref{q-wug},
we first establish several auxiliary lemmas.
\begin{lemma}\label{yl1}
Let $Q\subset\mathbb{R}^n$ be a cube, and let $k\in\mathbb{Z}_+$.
Then there exists a positive constant $C_{n,k}$ such that,
for any polynomial $P\in\mathcal{P}_k$,
\begin{equation}\label{yl1equ}
\|P\|_{L^{\infty}(Q)}\leq C_{n,k}\fint_{Q}
\left|P(x)\right|\,dx.
\end{equation}
\end{lemma}
\begin{proof}
Let $Q_0:=[-\frac12,\frac12]^n$.
Since $\mathcal{P}_{k}$ is finite-dimensional,
we infer that there exists an orthonormal basis
$\{\phi_i\}_{i=1}^N$ of $\mathcal{P}_{k}$ with respect
to the inner product
$$\left(R,S\right):=\int_{Q_0}R(x)\overline{S(x)}\,dx.$$
Hence, for any $P\in\mathcal{P}_{k}$ and $x\in Q_0$,
$$P(x)=\sum_{i=1}^N\int_{Q_0}P(y)\overline{\phi_i(y)}\,dy\phi_i(x).$$
Consequently,
\begin{align*}
|P(x)|&\leq \sum_{i=1}^N
\left|\int_{Q_0}P(y)\overline{\phi_i(y)}\,dy\right||\phi_i(x)|
\leq \sum_{i=1}^N
\|\phi_i\|_{L^\infty(Q_0)}
\int_{Q_0}|P(y)||\phi_i(y)|\,dy\\
&\leq\left[\sum_{i=1}^N\|\phi_i\|_{L^\infty(Q_0)}^2\right]
\int_{Q_0}|P(y)|\,dy.
\end{align*}
Therefore,
\begin{equation}\label{q0gj}
\|P\|_{L^\infty(Q_0)}
\leq C_{n,k}\int_{Q_0}|P(y)|\,dy
= C_{n,k}\fint_{Q_0}|P(y)|\,dy.	
\end{equation}

Now, let $Q$ be an arbitrary cube with center
$y_Q$ and edge length $\ell(Q)$. For any $z\in\mathbb{R}^n$,
let $\widetilde P(z):=P(y_Q+\ell(Q)z)$.
Then $\widetilde P\in\mathcal{P}_{k}$ and
$\|\widetilde P\|_{L^\infty(Q_0)}=\|P\|_{L^\infty(Q)}.$
Moreover, the change of variables $x=y_Q+\ell(Q)z$ gives
$$\fint_{Q_0}\left|\widetilde P(z)\right|\,dz=\fint_Q|P(x)|\,dx.$$
Applying \eqref{q0gj} with $P$ replaced by $\widetilde{P}$,
we obtain
$$\|P\|_{L^\infty(Q)}
=\left\|\widetilde P\right\|_{L^\infty(Q_0)}
\leq C_{n,k}\fint_{Q_0}\left|\widetilde P(z)\right|\,dz
=C_{n,k}\fint_Q|P(x)|\,dx.$$
This completes the proof of Lemma \ref{yl1}.
\end{proof}
The next lemma gives a version of the
Lebesgue differentiation
theorem for minimizing polynomials (see, for instance,
\cite[Lemma 4.1]{ds}).
\begin{lemma}\label{yl2}
Let $Q_0\subset\mathbb{R}^n$ be a
cube and $f\in L^1_{\rm{loc}}$. For any $x\in Q_0$ and
any $m\in\mathbb{N}$, let $Q_m^x$ denote
a dyadic subcube of
$Q_0$ with side length
$2^{-m}\ell(Q_0)$ containing $x$.
Then, for any $k\in\mathbb{Z}_+$
and almost every $x\in Q_0$,
$$\lim_{m\to\infty}P_{Q_m^x}^{(k)}f(x)=f(x),$$
where $P_Q^{(k)}f$ is the k-th
minimizing polynomial of $f$ on $Q$.
\end{lemma}
We next estimate the difference between
minimizing polynomials associated with two nested cubes.
\begin{lemma}\label{yl3}
Let $k\in\mathbb{Z}_+$, $Q\subset\mathbb{R}^n$ be a cube,
and $f\in L^1_{\rm{loc}}$.
If $Q'\subset Q$ is a cube satisfying $2\ell(Q')=\ell(Q)$, then
there exists a positive constant $C_{n,k}$, independent
of $f,$ $Q,$ and $Q'$, such that
\begin{equation}\label{zftjxdxs}
\left\|P^{(k)}_{Q'}(f)-P^{(k)}_{Q}(f)\right\|_{L^\infty(Q')}
\leq C_{n,k}\inf_{P\in\mathcal{P}_{k}}\fint_{Q}\left|
f(x)-P(x)\right|\,dx.	
\end{equation}
\end{lemma}
\begin{proof}
By Lemma \ref{jxyzx}, we conclude that
$$\fint_{Q'}\left|f(x)-P_{Q}^{(k)}(f)(x)\right|\,dx\leq
\frac{|Q|}{|Q'|}\fint_{Q}\left|f(x)-P_{Q}^{(k)}(f)(x)\right|
\,dx\lesssim
\inf_{P\in\mathcal{P}_{k}}\fint_{Q}\left|
f(x)-P(x)\right|\,dx.$$
From this and \eqref{yl1equ}, we deduce that
\begin{align*}
\left\|P^{(k)}_{Q'}(f)-P^{(k)}_{Q}(f)\right\|_{L^\infty(Q')}
&\lesssim\fint_{Q'}\left|P^{(k)}_{Q'}(f)(x)
-P^{(k)}_{Q}(f)(x)\right|\,dx\\
&\lesssim
\fint_{Q'}\left|P^{(k)}_{Q'}(f)(x)-f(x)\right|\,dx+
\fint_{Q'}\left|f(x)
-P^{(k)}_{Q}(f)(x)\right|\,dx\\
&\lesssim\fint_{Q'}\left|f(x)-P_{Q}^{(k)}(f)(x)\right|\,dx
\lesssim\inf_{P\in\mathcal{P}_{k}}\fint_{Q}\left|
f(x)-P(x)\right|\,dx,
\end{align*}
which completes the proof of Lemma \ref{yl3}.
\end{proof}
Combining the preceding lemmas, we obtain the following
pointwise estimate.
\begin{lemma}\label{yl4}
Let $k\in\mathbb{Z}_+$, $Q\subset\mathbb{R}^n$ be a cube,
and $f\in L^1_{\rm{loc}}$.
Then, for almost every $x\in Q$,
\begin{equation*}
\left|f(x)-P^{(k)}_{Q}(f)(x)\right|\leq C_{n,k}
\sum_{m=0}^{\infty}\inf_{P\in\mathcal{P}_k}\fint_{Q^x_{m}}
\left|f(y)-P(y)\right|\,dy,
\end{equation*}
where $Q_m^x$ denotes a dyadic subcube of $Q$ with edge length
$2^{-m}\ell(Q)$ containing $x$ and where $Q_0^x:=Q$.
\end{lemma}
\begin{proof}
By Lemma \ref{yl2}, for almost every $x\in Q$, we have
$$\left|f(x)-P_Q^{(k)}(f)(x)\right|=
\left|\sum_{m=0}^{\infty}\left[
P_{Q_{m+1}^x}^{(k)}(f)(x)-P_{Q_m^x}^{(k)}(f)(x)\right]\right|
\leq\sum_{m=0}^{\infty}\left|
P_{Q_{m+1}^x}^{(k)}(f)(x)-P_{Q_m^x}^{(k)}(f)(x)\right|.$$
Combining this with \eqref{zftjxdxs}, we obtain
$$\left|f(x)-P_Q^{(k)}(f)(x)\right|\lesssim
\sum_{m=0}^{\infty}\inf_{P\in\mathcal{P}_k}
\fint_{Q_m^x}\left|f(y)-P(y)\right|\,dy.$$
This completes the proof of Lemma \ref{yl4}.
\end{proof}
\begin{remark}
A corresponding result to Lemma \ref{yl4}
for balls can be found in \cite[Lemma 3.10]{hlyyz}.
\end{remark}
We are now ready to prove Theorem \ref{q-wug}.
\begin{proof}[Proof of Theorem \ref{q-wug}]
We first show the sufficiency (``if'' part). Let $\alpha\in(0,\infty)$.
By H\"older's inequality, we easily conclude that,
for any $q\in[1,p]$,
\begin{equation}\label{jnbz}
JN_{(p,p,k-1)_\alpha}^{\rm{con}}
\subset JN_{(p,q,k-1)_\alpha}^{\rm{con}}
\subset JN_{(p,1,k-1)_\alpha}^{\rm{con}}.	
\end{equation}
Thus, it remains
to show that $JN_{(p,1,k-1)_\alpha}^{\rm{con}}\subset
JN_{(p,p,k-1)_\alpha}^{\rm{con}}$.
Fix $f\in JN_{(p,1,k-1)_\alpha}^{\rm{con}}$, and
let $\{Q^j\}_j$
be an arbitrary collection of pairwise disjoint cubes with
common edge length $\ell$. Let $R_m$ denote the family
of all cubes obtained by subdividing each $Q^j$ dyadically
$m$ times. For any $x\in\bigcup_{j}Q^j$ and
$m\in\mathbb{N}$, there exists a unique cube $Q_m^x\in R_m$
such that $x\in Q_m^x$. Moreover, there exists a unique cube $Q^j$
such that $Q_m^x\subset Q^j$. From Lemma \ref{yl4} and
the triangle inequality for the norm $\{\sum_{j}\|\cdot\|^p_{L^p(Q_j)}\}^{\frac{1}{p}}$, we deduce that
\begin{align*}
&\left[\sum_{j}\int_{Q^j}\left|f(x)-P_{Q^j}^{(k-1)}(f)(x)
\right|^p
\,dx\right]^{\frac{1}{p}}\\
&\quad\lesssim
\left\{\sum_{j}\int_{Q^j}\left[\sum_{m=0}^{\infty}
\inf_{P\in\mathcal{P}_{k-1}}
\fint_{Q_{m}^x}\left|f(y)-P(y)\right|\,dy\right]^p
\,dx\right\}^{\frac{1}{p}}\\
&\quad\lesssim
\sum_{m=0}^{\infty}\left\{\sum_{j}\int_{Q^j}
\left[\inf_{P\in\mathcal{P}_{k-1}}
\fint_{Q_{m}^x}\left|f(y)-P(y)\right|\,dy\right]^p
\,dx\right\}^{\frac{1}{p}}.
\end{align*}
For $m\in\mathbb{N}$, let
$$g_m(x):=\inf_{P\in\mathcal{P}_{k-1}}
\fint_{Q_{m}^x}\left|f(y)-P(y)\right|\,dy.$$
Observe that $g_m$ is constant
on each cube in $R_m$. Consequently,
\begin{align*}
\sum_{m=0}^{\infty}\left\{\sum_{j}\int_{Q^j}
\left[g_m(x)\right]^p
\,dx\right\}^{\frac{1}{p}}&=\sum_{m=0}^{\infty}
\left\{\sum_{Q'\in R_m}\int_{Q'}
\left[g_m(x)\right]^p\,dx\right\}^{\frac{1}{p}}\\
&=\sum_{m=0}^{\infty}\left\{\sum_{Q'\in R_m}
|Q'|\left[\inf_{P\in\mathcal{P}_{k-1}}
\fint_{Q'}\left|f(y)-P(y)\right|\,dy\right]^p\,\right\}
^{\frac{1}{p}}\\
&\lesssim\|f\|_{JN^{\rm{con}}_{(p,1,k-1)_{\alpha}}}l
^{\alpha n}\sum_{m=0}^{\infty}2^{-mn\alpha}
\sim\|f\|_{JN^{\rm{con}}_{(p,1,k-1)_{\alpha}}}l
^{\alpha n}.
\end{align*}
From this, it follows that
\begin{align*}
\left[\sum_{j}|Q^j|^{1-p\alpha}
\fint_{Q_j}\left|f(x)-P_{Q^j}^{(k-1)}(f)(x)\right|^p
\,dx\right]^{\frac{1}{p}}\lesssim
\|f\|_{JN^{\rm{con}}_{(p,1,k-1)_{\alpha}}}.
\end{align*}
Taking the supremum over $\{Q^j\}_j\in\Pi_l$ and $l\in(0,\infty)$,
we obtain
$$\|f\|_{JN_{(p,p,k-1)_\alpha}^{\rm{con}}}
\lesssim\|f\|_{JN_{(p,1,k-1)_\alpha}^{\rm{con}}}.$$
Therefore, $JN_{(p,1,k-1)_\alpha}^{\rm{con}}\subset
JN_{(p,p,k-1)_\alpha}^{\rm{con}},$
which, together with \eqref{jnbz}, further implies that
$JN_{(p,1,k-1)_\alpha}^{\rm{con}}=
JN_{(p,p,k-1)_\alpha}^{\rm{con}}$. This completes the proof of the sufficiency.

We next prove the necessity (``only if'' part).
Let $\alpha=0$, and fix $q\in(1,p)$.
Then $\frac{1}{p-1}<\frac{q}{p-q}$, so
we can choose $\beta\in(\frac{1}{p-1},\frac{q}{p-q}]$.
For any $j\in\mathbb{N}$,
define $t_j:=\lfloor\beta\text{log}_2(j+1)\rfloor$,
$M_j:=M_{j-1}+t_j+2$, and $M_0:=1$.
This implies that, for any $j_1,j_2\in\mathbb{N}$
and $j_1\neq j_2$,
$$[M_{j_1}-t_{j_1},M_{j_1}]\cap[M_{j_2}-t_{j_2},M_{j_2}]
=\emptyset.$$
Let $\delta_j:=2^{-M_j}$, $L_j:=2^{-M_j+t_j}$, and
$$E_j:=\bigcup_{v=0}^{L_{j}^{-1}-1}\left\{
[j+vL_j,j+vL_j+\delta_j)\times[0,1)^{n-1}
\right\}\subset[j,j+1)\times[0,1)^{n-1}=:Q_j.$$
Clearly, the sets in this union are disjoint. Thus,
for any $j\in\mathbb{N}$,
$|E_j|=L_j^{-1}\delta_j=2^{-t_j}$.
We define
$$\partial E_j:=\left\{x\in\overline{E_j}:
\forall r\in(0,\infty), B(x,r)\cap E_j^\complement\neq\emptyset\right\}.$$
Choose $j_0\in\mathbb{N}$ such that, for any $j\geq j_0$,
$L_j>2\delta_j$ and $\delta_{j+1}\leq\frac{1}{4}\delta_{j}$.
Let
\begin{align}\label{equ-JN-def-f}
	f:=\sum_{j=j_0}^{\infty}f_j:=
\sum_{j=j_0}^{\infty}2^{\frac{t_j}{p}}
\mathbf{1}_{E_j}.
\end{align}
We first show that
$f\in JN^{\rm{con}}_{(p,1,0)_0}$. To this end,
fix $l\in(0,\infty)$ and $\mathcal{Q}\in\Pi_l$, and let
$$\text{I}_l:=\sum_{Q\in\mathcal{Q}}|Q|\left[\fint_{Q}
\left|f(y)-f_{Q}\right|\,dy\right]^p.$$
To estimate $\text{I}_l$, we
consider the following two cases for $\ell$:

\emph{Case (1): $l\leq1$.}
In this case, each cube $Q\in\mathcal{Q}$ intersects at most
two of the unit cubes $\{Q_j\}_{j\in\mathbb{N}}$
in sets of positive measure. Consequently,
\begin{equation}\label{Il}
\text{I}_l\lesssim\sum_{j=j_0}^{\infty}
\sum_{Q\in\mathcal{Q}}|Q|\left[\fint_{Q}
\left|f_j(y)-{f_j}_{Q}\right|\,dy\right]^p=:\sum_{j=j_0}^{\infty}
\text{J}_j.
\end{equation}
We estimate these terms by considering the following three cases for $\ell$:

\, \emph{Case (1.1): $0<l\leq\delta_j$.}
In this case, $l\leq\delta_j<L_j-\delta_j$.
If $Q\cap\partial E_j
=\emptyset$, then $f_j$ is constant almost
everywhere on $Q$. This implies $|f_j-{f_j}_Q|_Q=0.$
Therefore,
\begin{align}\label{J_l}
\text{J}_j&=\sum_{\{Q\in\mathcal{Q},Q\cap\partial E_j\neq\emptyset\}}
|Q|\left[\fint_{Q}
\left|f_j(y)-{f_j}_{Q}\right|\,dy\right]^p\notag\\
&=\sum_{\{Q\in\mathcal{Q},Q\cap\partial E_j\neq\emptyset\}}
|Q|^{1-p}\left[\int_{Q\cap E_j}
\left|2^{\frac{t_j}{p}}-2^{\frac{t_j}{p}}\frac{|Q\cap E_j|}{|Q|}
\right|\,dy+\int_{Q\setminus E_j}
\left|2^{\frac{t_j}{p}}\frac{|Q\cap E_j|}{|Q|}
\right|\,dy\right]^p\notag\\
&\lesssim\sum_{\{Q\in\mathcal{Q},Q\cap\partial E_j\neq\emptyset\}}
|Q|^{1-p}2^{t_j}
\left[\frac{|Q\setminus E_j||Q\cap E_j|}{|Q|}\right]^p\lesssim
\sum_{\{Q\in\mathcal{Q},Q\cap\partial E_j\neq\emptyset\}}
|Q|2^{t_j}.
\end{align}
Every cube in the last sum is contained in the set
$\{y:d(\partial E_j,y)<2\mathop{\rm{diam}}Q=2\sqrt{n}l\}$. Thus,
$$\#\{Q\in\mathcal{Q}:Q\cap\partial E_j\neq\emptyset\}
\leq\frac{|\{y:d(\partial E_j,y)<2\sqrt{n}l\}|}{|Q|}
\lesssim\frac{(4\sqrt{n}l+1)^{n-1}l}{L_jl^n}\lesssim
\frac{1}{L_jl^{n-1}}.$$
From this and \eqref{J_l}, we
deduce that
\begin{equation}\label{Jj1}
\text{J}_j\lesssim\frac{2^{t_j}l^n}{L_jl^{n-1}}\sim\frac{l}{\delta_j}.
\end{equation}

\, \emph{Case (1.2): $\delta_j<l\leq L_j$.}
In this case, for any $Q\in\mathcal{Q}$, there exist at most two indices $v$ such that
$$Q\cap[j+vL_j,j+vL_j+\delta_j)\times[0,1)^{n-1}\neq\emptyset.$$
Hence,
$|Q\cap E_j|\leq2\delta_jl^{n-1}$. This implies that
\begin{equation}\label{1.2gj}
|Q|\left[\fint_{Q}\left|f_j(y)-{f_j}_Q\right|\,dy\right]^p
\lesssim|Q|\left({f_j}_Q\right)^p\sim
|Q|2^{t_j}\left(\frac{|Q\cap E_j|}{|Q|}\right)^p\lesssim
l^n2^{t_j}\frac{\delta_j^p}{l^p}.	
\end{equation}
For any cube $Q\in\mathcal{Q}$ satisfying $|Q\cap E_j|>0$, we have
$$Q\subset\bigcup_{v=0}^{L_{j}^{-1}-1}\left\{
[j+vL_j-l,j+vL_j+\delta_j+l)\times[-l,1+l)^{n-1}
\right\}=:F_j.$$
From this and $\delta_j\leq l\leq1$, it follows that
$$\#\{Q\in\mathcal{Q}:|Q\cap E_j|>0\}
\leq\frac{|F_j|}{|Q|}\leq\frac{(2l+\delta_j)(2l+1)^{n-1}}{L_jl^n}
\lesssim\frac{l}{L_jl^n}.$$
Combining this with \eqref{1.2gj}, we obtain
\begin{equation}\label{Jj2}
\text{J}_j\lesssim\#\{Q\in\mathcal{Q}:|Q\cap E_j|>0\}
l^n2^{t_j}\frac{\delta_j^p}{l^p}
\lesssim2^{t_j}\frac{l\delta_j^p}{L_{j}l^p}
\lesssim\frac{\delta_j^{p-1}}{l^{p-1}}.
\end{equation}

\, \emph{Case (1.3): $L_j<l\leq1$.}
In this case, for any $Q\in\mathcal{Q}$, we have
$$\#\left\{v:[j+vL_j,j+vL_j+\delta_j)\times[0,1)^{n-1}\cap Q
\neq\emptyset\right\}\leq \frac{l}{L_j}+2.$$
Consequently,
\begin{equation}\label{1.3gj}
|Q|\left[\fint_{Q}\left|f_j(y)-{f_j}_Q\right|\,dy\right]^p
\lesssim|Q|\left({f_j}_Q\right)^p\sim
|Q|2^{t_j}\left(\frac{|Q\cap E_j|}{|Q|}\right)^p\lesssim
l^n2^{t_j}\frac{\delta_j^p(\frac{l}{L_j}+2)^p}{l^p}
\lesssim l^n2^{t_j}\frac{\delta_j^p}{L_j^p}.	
\end{equation}
For this estimate, it suffices to consider those cubes
$Q\in\mathcal{Q}$ satisfying
$|Q\cap E_j|>0$. Every such cube is contained in
$$[j-l,j+1+l]\times[-l,1+l]^{n-1}.$$
From this, we deduce that
$$\#\{Q\in\mathcal{Q}:|Q\cap E_j|>0\}\leq\frac{(2l+1)^n}{l^n}.$$
By this, $l\leq1$, and \eqref{1.3gj}, we find that
\begin{equation}\label{Jj3}
\text{J}_j\lesssim\frac{(2l+1)^n}{l^n}
l^n2^{t_j}\frac{\delta_j^p}{L_j^p}
\sim2^{t_j}\frac{\delta_j^p}{L_j^p}
\sim\frac{\delta_j^{p-1}}{L_j^{p-1}}.
\end{equation}
From \eqref{Jj1}, \eqref{Jj2}, \eqref{Jj3},
and \eqref{Il}, it follows that
\begin{equation}\label{Ilgj}
\text{I}_l\lesssim\sum_{\{j\geq j_0,l\leq\delta_j\}}\frac{l}{\delta_j}
+\sum_{\{j\geq j_0,\delta_j<l\leq L_j\}}\frac{\delta_j^{p-1}}{l^{p-1}}
+\sum_{\{j\geq j_0,L_j<l\}}\frac{\delta_j^{p-1}}{L_j^{p-1}}=:\text{S}_1
+\text{S}_2+\text{S}_3.
\end{equation}
We first estimate $\text{S}_1$. Let $G_l:=\{j\ge j_0:\delta_j\geq l\}$,
which is finite because $\delta_j\to0$ as $j\to\infty$.
If $G_l=\emptyset$, then $\text{S}_1=0$.
Otherwise, let $s_l:=\max G_l$. From
$\delta_{j+1}\leq\frac{1}{4}\delta_{j}$ and $\delta_j\geq l$,
we infer that
\begin{equation}\label{S1}
\text{S}_1\leq l\left(\frac{1}{\delta_{j_0}}
+...+\frac{1}{\delta_{s_l}}\right)
\leq\frac{l}{\delta_{s_l}}
\left(\frac{1}{4^{s_l}}+...+1\right)
\leq\sum_{j=0}^{\infty}\frac{1}{4^j}<\infty.
\end{equation}
Since the intervals
$\{[\delta_{j},L_{j}]\}_{j\in\mathbb{N}}$ are pairwise disjoint,
$\text{S}_2$ contains at most one term.
Using this and $\delta_j<l$, we obtain
\begin{equation}\label{S2}
\text{S}_2\lesssim1.	
\end{equation}
Finally, since $\beta(p-1)>1$, it follows that
\begin{equation}\label{S3}
\text{S}_3\lesssim\sum_{j=0}^{\infty}(j+1)^{-\beta(p-1)}<\infty.
\end{equation}
Combining \eqref{Ilgj}, \eqref{S1}, \eqref{S2},
and \eqref{S3}, we find that
\begin{equation}\label{Ilsj1}
\sup_{l\in(0,1]}\text{I}_l<\infty.
\end{equation}

\emph{Case (2): $l>1$.}
In this case, for any $Q\in\mathcal{Q}$, we have
$$\#\{j:Q\cap Q_j\neq\emptyset\}\leq\lfloor l\rfloor+2\leq3l.$$
Using this and \eqref{basic-inequality},
we conclude that
\begin{align*}
\text{I}_l&=\sum_{Q\in\mathcal{Q}}|Q|\left[\fint_Q|f(y)-f_Q|\,dy\right]^p
\lesssim \sum_{Q\in\mathcal{Q}}l^{n-np}\left[\int_Qf(y)\,dy\right]^p\\
&\lesssim\sum_{Q\in\mathcal{Q}}l^{n-np}\left(\sum_{j=j_0}^{\infty}
2^{\frac{t_j}{p}}|Q\cap E_j|\right)^p\\
&\lesssim l^{n-np}
\sum_{Q\in\mathcal{Q}}
(\#\{j:Q\cap Q_j\neq\emptyset\})^{p-1}\sum_{j=j_0}^{\infty}
2^{t_j}|Q\cap E_j|^p
\lesssim l^{(1-n)(p-1)}
\sum_{j=j_0}^{\infty}\sum_{Q\in\mathcal{Q}}2^{t_j}
|Q\cap E_j|^p.
\end{align*}
From this, $l>1$, and \eqref{basic-inequality} with
$\frac{1}{p}\in(0,1]$, we deduce that
\begin{align*}
\text{I}_l&\lesssim l^{(1-n)(p-1)}
\sum_{j=j_0}^{\infty}\sum_{Q\in\mathcal{Q}}2^{t_j}
|Q\cap E_j|^p\lesssim\sum_{j=j_0}^{\infty}2^{t_j}
\left(\sum_{Q\in\mathcal{Q}}|Q\cap E_j|\right)^p\\
&\lesssim\sum_{j=j_0}^{\infty}2^{t_j}|E_j|^p
\lesssim\sum_{j=j_0}^{\infty}2^{t_j(1-p)}
\lesssim\sum_{j=0}^{\infty}(j+1)^{-\beta(p-1)}<\infty.
\end{align*}
This implies that
\begin{equation}\label{Ilsj2}
\sup_{l\in(1,\infty)}\text{I}_l<\infty.
\end{equation}
By \eqref{Ilsj1} and \eqref{Ilsj2}, we obtain
$$\|f\|_{JN^{\rm{con}}_{(p,1,0)_0}}<\infty.$$
Hence, $f\in JN^{\rm{con}}_{(p,1,0)_0}$.
Using Lemma \ref{jxyzx}, we find that
$$\|f\|_{JN^{\rm{con}}_{(p,1,k-1)_0}}
\lesssim\|f\|_{JN^{\rm{con}}_{(p,1,0)_0}}<\infty.$$
Therefore,
$f\in JN^{\rm{con}}_{(p,1,0)_0}\subset
JN^{\rm{con}}_{(p,1,k-1)_0}$.

Finally, we prove that $f\notin
JN^{\rm{con}}_{(p,q,k-1)_0}$. From \eqref{jxdxsjfgj},
we infer that
$$\left\|P_{Q_j}^{(k-1)}(f)\right\|
_{L^{\infty}(Q_j)}\leq C_{n,k-1}\fint_{Q_j}|f(y)|\,dy
=C_{n,k-1}2^{\frac{t_j}{p}}\frac{|E_j|}{|Q_j|}=
C_{n,k-1}2^{t_j(\frac{1}{p}-1)}.$$
Consequently, for any $x\in E_j$,
\begin{align*}
\left|f(x)-P_{Q_j}^{(k-1)}(f)(x)\right|
\geq|f(x)|-\left\|P_{Q_j}^{(k-1)}(f)\right\|_{L^{\infty}(Q_j)}
\geq2^{\frac{t_j}{p}}\left(1-C_{n,k-1}2^{-t_j}\right).
\end{align*}
By this and the fact that $1-C_{n,k-1}2^{-t_j}\geq\frac{1}{2}$
whenever $j\ge\max\{\lfloor(4C_{n,k-1})^{\frac{1}{\beta}}\rfloor
,j_0\}=:j_1$ and $\beta\frac{p-q}{q}\leq1$, we obtain
\begin{align*}
\|f\|_{JN^{\rm{con}}_{(p,q,k-1)_0}}^p
&\geq\sum_{j=j_1}^{\infty}
|Q_j|\left[\fint_{Q_j}\left|f(x)-P_{Q_j}^{k-1}(f)(x)
\right|^q\,dx\right]^{\frac{p}{q}}\\
&\geq\sum_{j=j_1}^{\infty}
|E_j|^{\frac{p}{q}}2^{t_j}
\left(1-C_{n,k-1}2^{-t_j}\right)^p\\
&\gtrsim\sum_{j=j_1}^{\infty}2^{t_j(1-\frac{p}{q})}
\sim\sum_{j=j_1}^{\infty}(j+1)^{\beta(\frac{q-p}{q})}=\infty.
\end{align*}
This implies that
$$JN^{\rm{con}}_{(p,p,k-1)_0}\subset JN^{\rm{con}}_{(p,q,k-1)_0}
\subsetneqq JN^{\rm{con}}_{(p,1,k-1)_0},$$
which completes the proof of Theorem \ref{q-wug}.
\end{proof}

\noindent\textbf{Acknowledgements}\quad The main mathematical 
ideas and approaches in this article
originated from the authors. Generative AI tools were used as
supplementary aids during the exploratory and manuscript-preparation
stages, particularly in refining the presentation and examining
technical details. All mathematical arguments and proofs in the
final manuscript were independently verified and finalized by
the authors, who take full responsibility for their correctness.

\medskip

\noindent\textbf{Funding}\quad 
This project is partially supported by the National Natural Science
Foundation of China (Grant Nos. 12431006,
12371093, and 12501118), the Beijing Natural
Science Foundation (Grant No. 1262011), and the Fundamental Research
Funds for the Central Universities (Grant No. 2253200028).

\medskip

\noindent\textbf{Data Availability}\quad There is no data associated 
with the research in this article.

\section*{Declarations}

\noindent\textbf{Competing Interests}\quad The authors declare 
that they have no competing interests related to the research in
this article.

\bigskip

\smallskip

\noindent Zhexu Bai, Hao Peng, Dachun Yang,
Wen Yuan and Yangyang Zhang

\smallskip

\noindent Laboratory of Mathematics and Complex Systems (Ministry of Education of China),
School of Mathematical Sciences, Institute for Advanced Study, Beijing Normal University, Beijing 100875, The People's Republic of China

\smallskip

\noindent{\it E-mails}: \texttt{zhexubai@mail.bnu.edu.cn} (Z. Bai)

\noindent\phantom{\it E-mails:}  \texttt{haop@mail.bnu.edu.cn} (H. Peng)

\noindent\phantom{\it E-mails:} \texttt{dcyang@bnu.edu.cn} (D. Yang)

\noindent\phantom{\it E-mails:} \texttt{wenyuan@bnu.edu.cn} (W. Yuan)

\noindent\phantom{\it E-mails:} \texttt{yangyzhang@bnu.edu.cn} (Y. Zhang)

\end{document}